\documentclass{amsart}

\usepackage[latin1]{inputenc}
\usepackage[T1]{fontenc}
\usepackage{amsfonts}
\usepackage{amsmath}
\usepackage{amssymb}
\usepackage{graphicx}
\usepackage{pstricks}
\usepackage{pst-grad}
\usepackage{multido}
\usepackage{amsthm}

\newcommand{\cA}{\mathcal{A}}
\newcommand{\cF}{\mathcal{F}}

\newcommand{\bZ}{\mathbb{Z}}
\newcommand{\bN}{\mathbb{N}}
\newcommand{\sg}{\mathcal{h}}
\newcommand{\sd}{\mathcal{i}}
\newcommand{\la}{\lambda}
\newcommand{\cH}{\mathcal{H}}

\numberwithin{equation}{section}
\newtheorem{Th}{Theorem}[section]
\newtheorem{Cor}[Th]{Corollary}
\newtheorem{Lem}[Th]{Lemma}
\newtheorem{Def}[Th]{Definition}

\newtheorem{Prop}[Th]{Proposition}
\newtheorem{Rem}[Th]{Remark}
\newtheorem*{thm1}{Theorem 1.2}

\begin{document}

\title{On the determination of Kazhdan-Lusztig cells for affine Weyl groups with unequal parameters}
\author{J\'er\'emie Guilhot}
\address{Department of Mathematical Sciences, King's College,
Aberdeen University, Aberdeen AB24 3UE, Scotland, U.K.}
\address{Institut Camille Jordan \\Universit\'e Claude Bernard Lyon 1,Lyon,  France.}

\email{guilhot@maths.abdn.ac.uk}

\date{July, 2007}

\dedicatory{Dedicated to the memory of Fokko du Cloux}

\begin{abstract}
Let $W$ be a Coxeter group and $L$ be a weight function on $W$. Following
Lusztig, we have a corresponding decomposition of $W$ into left cells which have important applications in representation theory. We study
the case where $W$ is an affine Weyl group of type $\tilde{G_{2}}$. Using
explicit computation with \textsf{COXETER} and \textsf{CHEVIE} , we show that (1) there are
only finitely many possible decompositions into left cells and (2) the
number of left cells is finite in each case, thus confirming some of
Lusztig's conjectures in this case. A key ingredient of the proof is
a general result which shows that the Kazhdan-Lusztig polynomials of affine
Weyl group are invariant under (large enough) translations.  
\end{abstract}

\maketitle

\pagestyle{myheadings}
\markboth{Guilhot}{Kazhdan-Lusztig cells for affine Weyl groups with unequal parameters}

\maketitle


\section{Introduction}
This paper is concerned with Kazhdan-Lusztig polynomials and left cells in an affine Weyl group $W$ (with set of simple reflections $S$) with respect to a weight function $L$ where, following Lusztig \cite{bible}, a weight function on $W$ is a function $L:W \rightarrow \mathbb{Z}$ such that $L(ww')=L(w)+L(w')$ whenever $\ell(ww')=\ell(w)+\ell(w')$  ($\ell$ is the usual length function on $W$). We shall only consider weight functions such that $L(w)>0$ for all $w\neq 1$.

Generalized left, right and two-sided cells of a Coxeter group give rise to left, right and two-sided modules of the corresponding Hecke algebra $\mathcal{H}$ with parameters given by $L$. In turn, the representation theory of the Hecke algebra is very relevant to the representation theory of reductive groups over $p$-adic fields.

The case where $L$ is constant on the generators is known as the equal parameter case. It was studied by Lusztig in \cite{Lus3,Lus4,Lus5,Lus6}. The left cells have been explicitly described for type $\tilde{A}_{r},r\in\mathbb{N}$ (see \cite{Lus2, Shi1}), ranks 2, 3 (see \cite{Lus3, Bed,  Du}) and types $\tilde{B}_{4}$, $\tilde{C}_{4}$ and $\tilde{D}_{4}$ (see \cite{Shi2,Shi3,Chen}).

Much less is known for unequal parameters. In \cite{bible}, Lusztig has formulated a number of precise conjectures in that case. The proof of these conjectures in the equal parameter case involved an interpretation of the Kazhdan-Lusztig polynomials in terms of intersection cohomology which ensures that all the coefficients of these polynomials are non-negative integers. In the general case, this does not hold anymore. 

The case where the parameters are coming from a graph automorphism has been studied by K. Bremke in \cite{Bremke}, using the previous properties. 

One of the many consequences of Lusztig's conjectures is that the number of left cells is finite in an affine Weyl group. One of our aim here is to show the following theorem.
\begin{Th}
\label{finite}
Let $W$ be an affine Weyl group of type $\tilde{G}_{2}$ and $L$ be any weight function on $W$ with $L(w)>0$ for $w\neq 1$. Then the following hold.
\begin{enumerate}
\item There are
only finitely many possible decompositions of $W$ into left cells.
\item The
number of left cells is finite in each case.
\end{enumerate}
\end{Th}

For the proof we must show a sufficient number of equalities among Kazhdan-Lusztig polynomials (see Section 3 for precise definitions). These are provided by the following theorem which holds for any affine
Weyl group $W$. Let $u\in W$ be a translation (see Definition \ref{deftrans}). For $x,y\in W$ we write $x.y$ if and only if $\ell(xy)=\ell(x)+\ell(y)$. We have
\begin{Th}
\label{translation}
Let $z,z'\in W$ and $N=\ell(z')-\ell(z)\geq 0$. Then for $r> N(\ell(u)+1)$ and for any $k\geq 0$, we have
$$\mathbf{ P}_{z.u^{r},z'.u^{r}}=\mathbf{ P}_{z.u^{r+k},z'.u^{r+k}}$$
and if there exists $s\in S$ which satisfies $sz.u^{r}<z.u^{r}<z'.u^{r}<sz'.u^{r}$, then 
$$\mathbf{ M}_{z.u^{r},z'.u^{r}}=\mathbf{ M}_{z.u^{r+k},z'.u^{r+k}}\ .$$
\end{Th}

While the statement and the proof of Theorem 1.2 do not make any reference to computer calculations, we would like to point out that
we arrived at the precise statement through extensive experimentation
using Ducloux's COXETER program \cite{Cox}.

\section{A geometric realization of affine Weyl groups}
In this section, we recall some basic material about affine Weyl groups which will be needed later on. The exposition follows \cite{Lus1,Bremke,Xi} and we refer to these publications for more details and proofs. Some of the notions are illustrated in Figure 1, at the end of the section.

Let $V$ be an Euclidean space of finite dimension $r\geq 1$.
Let $\Phi\subset V$ be an irreducible root system of rank $r$ and $\check{\Phi}\subset V^{*}$ the dual root system. Fix a set of positive roots $\Phi^{+}\subset\Phi$. We denote the coroot corresponding to $\alpha$ by $\check{\alpha}$ and we write $\mathcal{h}x,y\mathcal{i}$ for the value of $y\in V^{*}$ at $x\in V$. Let $W_{0}$ be the Weyl group of $\Phi$. For $\alpha\in\Phi$ and $k\in \mathbb{Z}$, we define a hyperplane
$$H_{\alpha,k}=\{\lambda\in V\mid \mathcal{h}\lambda,\check{\alpha}\mathcal{i}=k\}.$$
Let 
$$\mathcal{F}=\{H_{\alpha,k}\mid \alpha\in \Phi^{+}, k\in \mathbb{Z}\}.$$ 
Each $H\in \mathcal{F}$ defines an orthogonal reflection $\sigma_{H}$ in $V$ with fixed point set $H$. Let $\Omega$ be the group of affine transformations generated by these reflections. We regard $\Omega$ as acting on the right on $V$. 
An alcove is a connected component of the set 
$$V-\underset{H\in\mathcal{F}}{\cup} H.$$
$\Omega$ acts simply transitively on the set $X$ of alcoves. 

Let $S$ be the set of $\Omega$-orbits in the set of faces (codimension 1 facets) of alcoves. Then $S$ consists of $r+1$ elements which can be represented as the $r+1$ faces of any given alcove.

For $s\in S$ we define an involution $A\rightarrow sA$ of $X$ as follows. Given an alcove A, we denote by $sA$ the unique alcove distinct from A which shares with $A$ a face of type $s$. The maps $A \rightarrow sA$ generate a group of permutations of $X$ which is a Coxeter group $(W,S)$. In our case, it is the affine Weyl group usually denoted by $\tilde{W_{0}}$. We regard $W$ as acting on the left on $X$. It acts simply transitively and it commutes with the action of $\Omega$ on $X$. 

We fix parameters $c_{s}\in\mathbb{N}^{*}$ for $s\in S$. Recall that, for $s,t\in S$, the parameters must only satisfy $c_{s}=c_{t}$ if $s$ and $t$ are conjugate in $W$. In the case where $W$ is of type $\tilde{C}_{r+1}$ with generators $s_{1},\ldots,s_{r+1}$ and $W_{0}$ is generated by $s_{1},\ldots,s_{r}$ we can assume, by the symmetry of the Dynkin diagram, that $c_{s_{1}}\geq c_{s_{r+1}}$. Similarly, if $W$ is of type $\tilde{A}_{1}$ with generators $s_{1},s_{2}$ and $W_{0}$ is generated by $s_{1}$, we can assume that $c_{s_{1}}\geq c_{s_{2}}$. In \cite{Bremke}, Bremke showed that if a hyperplane $H\in \mathcal{F}$ supports faces of types $s,t\in S$ then $s$ and $t$ are conjugate in $W$. 
As a consequence of this result, we can associate an integer $c_{H}\in\mathbb{N}$ to $H\in\mathcal{F}$, where $c_{H}=c_{s}$ if $H$ supports a face of type $s$.

For a $0$-dimensional facet $v$ of an alcove, let
$$m(v)=\underset{H\ni v,H\in\mathcal{F}}{\sum}c_{H}.$$
We say that $v$ is a special point if $m(v)$ is maximal. If $c_{s}=1$ for all $s\in S$, then the notion of special points is the same as the notion in \cite{Lus1} and $m(v)=|\Phi^{+}|$ for any special point. Note that, following \cite[Section 2]{Bremke}, and with our convention for $\tilde{C}_{r+1}$ and $\tilde{A}_{1}$, $0\in V$ is a special point.

Let $n\in\bZ$. A hyperplane $H=H_{\alpha,n}\in\mathcal{F}$ divides $V-H$ into two half-spaces:
$$\{x\in V\mid \mathcal{h}x,\check{\alpha}\mathcal{i}>n\}$$
and
$$\{x\in V\mid \mathcal{h}x,\check{\alpha}\mathcal{i}<n\}.$$
Let $v$ be a special point. A quarter with vertex $v$ is a connected component of the set
$$V-\underset{H\ni v,H\in\mathcal{F}}{\bigcup}H.$$
Hyperplanes which are adjacent to a quarter $\mathcal{C}$ are called walls of $\mathcal{C}$. 

Let $A_{0}$ be the fundamental alcove defined by
$$\{v\in V\mid 0<\mathcal{h}x,\check{\alpha}\mathcal{i}<1 \text{ for every positive root $\alpha$}\}.$$

If $\lambda$ is a $0$-dimensional facet of an alcove, we denote by $\Omega_{\lambda}$ the stabilizer of $\la$ in $\Omega$ and by $W_{\lambda}$ the stabilizer in $W$ of the set of alcoves $A$ which contain $\lambda$ in their closure. $W_{\lambda}$ is generated by $r$ elements of $S$. It is a maximal parabolic subgroup of $W$ (if $\lambda=0\in V$, the definition of $W_{\lambda}$ is consistent with the definition of $W_{0}$ given before).

We now introduce a new definition.
\begin{Def}
\label{hdef}
Let $z\in W$ and $A\in X$. Let $H_{1},...,H_{n}$ be the set of hyperplanes which separate $A$ and $zA$. For $1\leq i\leq n$, let $E_{H_{i}}(zA)$ be the half-space defined by $H_{i}$ which contains $zA$. 
Let
$$h_{A}(z)=\overset{n}{\underset{i=1}{\bigcap}} E_{H_{i}}(zA).$$
\end{Def}

For $z,z'\in W$ we write $z.z'$ if and only if $\ell(zz')=\ell(z)+\ell(z')$. It is well known that for any $w\in W$ and any $A\in X$, $\ell(w)$ is the number of hyperplanes which separate $A$ and $wA$. Therefore one can see that $z'.z$ if and only if
$$\{H\mid H \text{ separates } A \text{ and } zA\} \cap  \{H\mid H \text{ separates } zA \text{ and } z'zA\}=\emptyset,$$
or in other words
\begin{Lem}
\label{hzone}
Let $z,z'\in W$ and $A\in X$. We have
$$z.z'\Leftrightarrow z(z'A)\subset h_{A}(z').$$ 
\end{Lem}

In Figure 1, we consider an affine Weyl group $W$ of type $\tilde{G}_{2}$ 
$$W:=\sg s_{1},s_{2},s_{3}\ |\ (s_{1}s_{2})^{6}=1,(s_{2}s_{3})^{3}=1,(s_{1}s_{3})^{2}=1\sd.$$
The thick arrows represent the set of positive roots $\Phi^{+}$, $\lambda$ is a special point, the gray area around $\la$ is the set of alcoves containing $\la$ in their closure, $zA_{0}$ is the image of the fundamental alcove $A_{0}$ under the action of $z=s_{3}s_{2}s_{1}s_{2}s_{1}s_{2}\in W$ and $\mathcal{C}$ is a quarter with vertex $0$. Note that the subgroup $W_{\lambda}$ is generated by $s_{1}$ and $s_{2}$.

\begin{center}
\begin{pspicture}(-3,4.37)(3,-4.37)
\psset{unit=2.5cm}

\pslinewidth=.3pt

\pspolygon[fillcolor=lightgray,fillstyle=solid](0.866,0.5)(0.866,1)(1.3,1.25)(1.732,1)(1.732,0.5)(1.3,0.25)
\rput(1.26,0.62){{\footnotesize $\lambda$}}

\pspolygon[fillcolor=lightgray,fillstyle=solid](0,0)(0,0.5)(0.2166,0.375)
\rput(0.085 ,0.35){{\footnotesize$A_{0}$}}
\psline[linewidth=0.5mm]{->}(0,0)(-1.3,0.75)
\psline[linewidth=0.5mm]{->}(0,0)(0.866,0)
\psline[linewidth=0.5mm]{->}(0,0)(-.433,0.75)
\psline[linewidth=0.5mm]{->}(0,0)(0,1.5)
\psline[linewidth=0.5mm]{->}(0,0)(0.433,0.75)
\psline[linewidth=0.5mm]{->}(0,0)(1.3,0.75)

\psline[linewidth=0.5mm](-1.732,-1)(-0.433,-0.25)
\psline[linewidth=0.5mm](-0.433,-0.25)(0,-0.5)
\psline[linewidth=0.5mm](0,-0.5)(0,-1.5)
\rput(-0.7,-1.1){$h_{A_{0}}(z)$}


\pspolygon[fillcolor=lightgray,fillstyle=solid](0,0)(1.732,-1)(1.732,-1.5)(0.866,-1.5)
\rput(0.8,-0.65){{\footnotesize $\mathcal{C}$}}
\pslinewidth=.3pt

\rput(-0.08,-0.08){{\footnotesize $0$}}

\pspolygon[fillstyle=solid,fillcolor=lightgray](-0.433,-0.75)(-0.433,-0.25)(-0.2166,-0.375)
\rput(-0.33,-0.4){{\scriptsize $zA_{0}$}}

\psline(-1.732,1.5)(1.732,1.5)
\psline(-1.732,0.75)(1.732,0.75)
\psline(-1.732,0)(1.732,0)

\psline(-1.732,-0.75)(1.732,-0.75)
\psline(-1.732,-1.5)(1.732,-1.5)


\multido{\n=-1.732+0.433}{9}{
\psline(\n,-1.5)(\n,1.5)}

\psline(-1.732,1)(-0.866,1.5)
\psline(-1.732,0.5)(0,1.5)
\psline(-1.732,0)(0.866,1.5)
\psline(-1.732,-0.5)(1.732,1.5)

\psline(-1.732,-1)(1.732,1)

\psline(1.732,-1)(0.866,-1.5)
\psline(1.732,-0.5)(0,-1.5)
\psline(1.732,0)(-0.866,-1.5)
\psline(1.732,0.5)(-1.732,-1.5)


\psline(1.732,1)(0.866,1.5)
\psline(1.732,0.5)(0,1.5)
\psline(1.732,0)(-0.866,1.5)
\psline(1.732,-0.5)(-1.732,1.5)

\psline(1.732,-1)(-1.732,1)

\psline(-1.732,-1)(-0.866,-1.5)
\psline(-1.732,-0.5)(0,-1.5)
\psline(-1.732,0)(0.866,-1.5)
\psline(-1.732,0.5)(1.732,-1.5)

\psline(-1.732,0)(-0.866,1.5)
\psline(-1.732,-1.5)(0,1.5)
\psline(-0.866,-1.5)(0.866,1.5)
\psline(0,-1.5)(1.732,1.5)
\psline(0.866,-1.5)(1.732,0)

\psline(1.732,0)(0.866,1.5)
\psline(1.732,-1.5)(0,1.5)
\psline(0.866,-1.5)(-0.866,1.5)
\psline(0,-1.5)(-1.732,1.5)
\psline(-0.866,-1.5)(-1.732,0)
\rput(0,-1.75){fig 1. $\tilde{G}_{2}$}
\end{pspicture}
\end{center}

\section{Total ordering and weight function}
The basic references for this section are \cite{Lus1p,bible,geck}. In \cite{bible}, Lusztig studies the left cells of a Coxeter group $W$ with respect to a weight function $L$ on $W$.  Considering a more abstract setting as defined by Lusztig in \cite{Lus1p}, where left cells are defined with respect to an abelian group and a total order on it, Geck \cite{geck} formulates some conditions for two weight functions to give rise to the same left cell decomposition (when $W$ is finite). In this section we will find some conditions for two weight functions to give rise to essentially the same Kazhdan-Lusztig polynomials on a given subset of $W$ (when $W$ is an affine Weyl group). 

We recall the basic setting for the definition of Kazhdan-Lusztig polynomials and left cells. Let $W$ be a Coxeter group with generating set $S$.  Let $\Gamma$ be an abelian group (written multiplicatively) and $\mathbf{ A}:=\mathbb{Z}[\Gamma]$ be the group algebra of $\Gamma$ over $\mathbb{Z}$. Let $\{v_{s}\mid s\in S\}$ be a subset of $\Gamma$ such that $v_{s}=v_{t}$ whenever $s,t\in S$ are conjugate in $W$. Then we can define the corresponding generic Iwahori-Hecke algebra $\mathbf{ H}$, with $\mathbf{ A}$-basis $\{\mathbf{ T}_{w}\mid w\in W\}$ and multiplication given by the rule
\begin{equation*}
\mathbf{ T}_{s}\mathbf{ T}_{w}=
\begin{cases}
\mathbf{ T}_{sw}, & \mbox{if } \ell(sw)>\ell(w),\\
\mathbf{ T}_{sw}+(v_{s}-v_{s}^{-1})\mathbf{ T}_{w}, &\mbox{if } \ell(sw)<\ell(w);
\end{cases}
\end{equation*}
where $\ell: W\rightarrow \mathbb{N}$ denotes the usual length function on $W$ with respect to S. 

Let $a\rightarrow \bar{a}$ be the involution of $\mathbb{Z}[\Gamma]$ defined by $\overline{g}=g^{-1}$ for $g\in \Gamma$. We can extend it to a map from $\mathbf{ H}$ to itself by
\begin{equation*}
\overline{\underset{w\in W}{\sum}a_{w}\mathbf{ T}_{w}}=\underset{w\in W}{\sum}\bar{a}_{w}\mathbf{ T}^{-1}_{w^{-1}} \qquad (a_{w}\in\bZ[\Gamma]).
\end{equation*}
Then $h\rightarrow \bar{h}$ is a ring involution.

For $w\in W$, define $\mathbf{ r}_{y,w}\in \mathbf{ A}$ by
\begin{equation*}
\overline{\mathbf{ T}}_{w}=\underset{y\in W}{\sum}\bar{\mathbf{ r}}_{y,w}\mathbf{ T}_{y}.
\end{equation*}
These $\mathbf{ r}$-polynomials satisfy $\mathbf{ r}_{y,w}=0$ unless $y\leq w$, $\mathbf{ r}_{y,y}=1$ and the following recursive formula, for $s\in S$ such that $sw<w$ (where $<$ denotes the Bruhat order on $W$)
\begin{equation*}
\label{rec rpol}
\mathbf{ r}_{y,w}=
\begin{cases}
\mathbf{ r}_{sy,sw},&\text{if}\ sy<y,\\
\mathbf{ r}_{sy,sw}+(v_{s}-v_{s}^{-1})\mathbf{ r}_{y,sw},&\text{if}\ sy>y.\\
\end{cases}
\end{equation*}

Choose a total ordering of $\Gamma$. This is specified by a multiplicatively closed subset $\Gamma_{+}\subset \Gamma -\{1\}$ such that $\Gamma=\Gamma_{+}\cup \{1\}\cup \Gamma_{-}$ (disjoint union) where $\Gamma_{-}=\{g^{-1}\mid g\in \Gamma_{+}\}$. Moreover, assume that
$$\{v_{s}\mid s\in S\}\subset \Gamma_{+}.$$
Given a total ordering as above, there exists a unique element $\mathbf{ C}_{w}\in \mathbf{ H}$ such that
$$\mathbf{ C}_{w}=\overline{\mathbf{ C}}_{w}\ \ \text{  and  }\ \ \ \ \mathbf{ C}_{w}=\mathbf{ T}_{w}+
\underset{y\in W, y<w}{\sum}\mathbf{ P}_{y,w}\mathbf{T}_{y},$$
where $\mathbf{ P}_{y,w}\in \mathbb{Z}[\Gamma_{-}]$ for $y<w$. In fact, the set $\{\mathbf{ C}_{w}\mid w\in W\}$ forms a basis of $\mathbf{ H}$ known as the Kazhdan-Lusztig basis. We set $\mathbf{ P}_{y,y}=1$ for any $y\in W$.

The following formula gives the relation between the Kazhdan-Lusztig polynomials $\mathbf{ P}$ and the $\mathbf{ r}$-polynomials
\begin{equation}
\label{P dependence}
\overline{\mathbf{ P}}_{x,w}=\underset{y;x\leq y\leq w}{\sum}\mathbf{ r}_{x,y}\mathbf{ P}_{y,w}.
\end{equation}
Note that this formula together with the condition $\mathbf{ P}_{y,w}\in \Gamma_{-}$ for $y<w$, uniquely defines the Kazhdan-Lusztig polynomials.


Let $w\in W$ and $s\in S$, we have the following multiplication formula
$$\mathbf{ C}_{s}\mathbf{ C}_{w}=
\begin{cases}
\mathbf{ C}_{sw}+\underset{z;\ sz<z<w}{\sum}\mathbf{ M}_{z,w}^{s}\mathbf{ C}_{z}, & \mbox{if }w<sw,\\ 
(v_{s}+v_{s}^{-1})\mathbf{ C}_{w}, & \mbox{if }sw<w; 
\end{cases}
$$
where $\mathbf{M}_{z,w}^{s}\in \mathbb{Z}[\Gamma]$ satisfy
\begin{equation}
 \overline{\mathbf{ M}_{y,w}^{s}}=\mathbf{ M}_{y,w}^{s},
 \end{equation}
 \begin{equation}
 \label{M dependence}
(\underset{z;y\leq z<w;sz<z}{\sum}\mathbf{ P}_{y,z}\mathbf{ M}_{z,w}^{s})-v_{s}\mathbf{ P}_{y,w} \in \mathbb{Z}[\Gamma_{-}].
\end{equation}

Given $y,w\in W$ and $s\in S$ we write $y\leftarrow_{L} w$ if $\mathbf{ C}_{y}$ appears with a non-zero coefficient in $\mathbf{ C}_{s}\mathbf{ C}_{w}$ for some $s\in S$.
The Kazhdan-Lusztig pre-order $\leq_{L}$ is the transitive closure of the above relation, i.e $y\leq_{L} w$ if there exists a sequence $y=y_{0},...,y_{n}=w$ in $W$ such that $y_{i}\leftarrow_{L}y_{i+1}$ for all $0\leq i\leq n-1$. The equivalence relation associated to $\leq_{L}$ will be denoted 
by $\sim_{L}$ and the corresponding equivalence classes are called left cells of $W$.
\begin{Rem}
\label{P,M dependence} 
Let $y\leq w\in W$ and $[y,w]=\{z\in W\ |\ y\leq z\leq w\}$. Looking at the relations \ref{P dependence}--\ref{M dependence}, one can see that the set of polynomials
$$\{\mathbf{P}_{x,z},\  \mathbf{M}_{x,z}^{s}\ |\ x,z\in [y,w]\}$$
is determined by the total ordering of $\Gamma$ and the set of polynomials
$$\{\mathbf{r}_{z_{1},z_{2}}\ |\ z_{1},z_{2}\in [y,w]\}.$$
\end{Rem}
We now specialize the above setting to the case where the parameters of the Iwahori-Hecke algebra  are given by a weight function. 
We only consider weight functions such that $L(s)>0$ for all $s\in S$. Let $\cA=\mathbb{Z}[v,v^{-1}]$ where $v$ is an indeterminate. We have a corresponding Iwahori-Hecke algebra $\cH$ with parameters $\{v^{L(s)}\mid s\in S\}$. As before, $\cH$ has an $\cA$-basis $\{ T_{w}\mid w\in W\}$ with multiplication given by the formula
\begin{equation*}
T_{s}T_{w}=
\begin{cases}
T_{sw}, &\text{if $\ell(sw)>\ell(w)$,}\\
T_{sw}+(v^{L(s)}-v^{-L(s)})T_{w}, &\text{if $\ell(sw)<\ell(w).$}
\end{cases}
\end{equation*}

Now consider the abelian group $\{v^{n}\mid n\in \mathbb{Z}\}$ with the total order specified by $\{v^{n}\mid n>0\}$. Thus as above, we can define the Kazhdan-Lusztig basis $\{C_{w}\mid w\in W\}$ of $\cH$. 
We obtain
\begin{enumerate}
\item a collection of polynomials $r_{y,w}\in \mathbb{Z}[v,v^{-1}]$,
\item a collection of polynomials $P_{y,w}\in v^{-1}\mathbb{Z}[v^{-1}]$ for all $y<w \in W$,
\item a collection of polynomials $M_{y,w}^{s}\in \mathbb{Z}[v,v^{-1}]$ where $sy<y<w<sw$.
\end{enumerate}

In \cite{geck}, Geck has established a link between these two situations, where you have an abelian group $\Gamma$ with a total order specified by $\Gamma_{+} \subset \Gamma$ and a choice of parameters $\{ v_{s}\mid s\in S \} \subset \Gamma_{+}$ on the one hand, and a weight function L on the other hand. Denote by $\mathbf{ r}_{y,w}$, $\mathbf{ P}_{y,w}$ and $\mathbf{ M}_{y,w}^{s}$ the polynomials in $\mathbb{Z}[\Gamma]$ arising in the first case and by $r_{y,w}$, $P_{y,w}$ and $M_{y,w}^{s}$ the polynomials in $\mathbb{Z}[v,v^{-1}]$ arising in the second case.

He defined two subsets $\Gamma_{+}^{a}(W),\Gamma_{+}^{b}(W)\subset \Gamma_{+}$ as follows. First, let $\Gamma_{+}^{a}(W)$ be the set of all elements $\gamma \in \Gamma_{+}$ such that $\gamma^{-1}$ occurs with a non zero coefficient in a polynomial $\mathbf{P}_{y,w}$ for some $y<w$ in $W$. Next for any $y,w$ in $W$ and $s\in S$ such that $\mathbf{M}_{y,w}^{s}\neq 0$, we write $\mathbf{M}_{y,w}^{s}=n_{1}\gamma_{1}+...+n_{r}\gamma_{r}$ where $0\neq n_{i}\in \mathbb{Z}$, $\gamma_{i}\in \Gamma$ and $\gamma_{i-1}^{-1}\gamma_{i}\in \Gamma_{+}$ for $2\leq i\leq r$. Let $\Gamma_{+}^{b}(W)$ be the set of all elements  $\gamma_{i-1}^{-1}\gamma_{i}\in \Gamma_{+}$ arising in this way, for any $y,w,s$ such that $\mathbf{ M}_{y,w}^{s}\neq 0$. Finally set $\Gamma_{+}(W)=\Gamma_{+}^{a}(W)\cup\Gamma_{+}^{b}(W)$. Then he proved that
\begin{Prop}
Assume that we have a ring homomorphism 
$$\sigma: \mathbb{Z}[\Gamma]\rightarrow \ \mathbb{Z}[v,v^{-1}],\ \ v_{s}\rightarrow v^{L(s)}$$
such that
\begin{equation*}
\sigma(\Gamma_{+}(W))\subseteq \{v^{n}\mid n>0\}.\tag{$\ast$}
\end{equation*} 
Then $\sigma(\mathbf{ P}_{y,w})=P_{y,w}$ for all $y<w$ in W and $\sigma(\mathbf{ M}_{y,w}^{s})=M_{y,w}^{s}$ for any $y,w\in W$ such that $sy<y<w<sw$. Furthermore, the relation $\leq_{L}$ on $W$ defined with respect to the weight function $L$ is the same as the one defined with respect to $\Gamma_{+}\subset \Gamma$.
\end{Prop}

In order to deal with affine Weyl groups (infinite groups), we need a refinement of the above result for Bruhat intervals.

Let $y,w\in W$, $s\in S$ and $I=[y,w]$. We now define three subsets $\Gamma_{+}^{a}(I)$, $\Gamma_{+}^{b,s}(I)$, $\Gamma_{+}^{c,s}\subset \Gamma_{+}$. First, let $\Gamma_{+}^{a}(I)$ be the set of all elements $\gamma \in \Gamma_{+}$ such that $\gamma^{-1}$ occurs with a non-zero coefficient in a polynomial $\mathbf{P}_{z_{1},z_{2}}$ for some $z_{1}<z_{2}$ in $I$. Next for any $z_{1},z_{2}$ in $I$ such that $\mathbf{M}_{z_{1},z_{2}}^{s}\neq 0$ we write $\mathbf{M}_{z_{1},z_{2}}^{s}=n_{1}\gamma_{1}+...+n_{r}\gamma_{r}$ where $0\neq n_{i}\in \mathbb{Z}$, $\gamma_{i}\in \Gamma$ and $\gamma_{i-1}^{-1}\gamma_{i}\in \Gamma_{+}$ for $2\leq i\leq r$. Let $\Gamma_{+}^{b,s}(I)$ be set of all elements  $\gamma_{i-1}^{-1}\gamma_{i}\in \Gamma_{+}$ arising in this way, for any $z_{1},z_{2}\in I$ such that $\mathbf{ M}_{z_{1},z_{2}}^{s}\neq 0$. Finally let $\Gamma_{+}^{c,s}$ be the set of all elements $\gamma \in \Gamma_{+}$ such that $\gamma^{-1}$ occurs with a non-zero coefficient in a polynomial of the form 
$$\underset{z;z_{1}\leq z<z_{2};sz<z}{\sum}\mathbf{ P}_{z_{1},z}\mathbf{ M}_{z,z_{2}}^{s}-v_{s}\mathbf{ P}_{z_{1},z_{2}}$$
where $z_{1},z_{2}\in I$ and $sz_{1}<z_{1}<z_{2}<sz_{2}$.
We set $\Gamma_{+}^{s}(I)=\Gamma_{+}^{a}(I)\cup\Gamma_{+}^{b,s}(I)\cup \Gamma_{+}^{c,s}$.
\begin{Prop}
\label{L}
Let $y,w\in W$, $s\in S$ and $I=[y,w]$. Assume that we have a ring homomorphism 
\begin{equation*}
\sigma : \mathbb{Z}[\Gamma]\rightarrow \ \mathbb{Z}[v,v^{-1}],\ \ v_{s}\rightarrow v^{L(s)} 
\end{equation*}
such that 
\begin{equation*}
\sigma(\Gamma_{+}^{s}(I))\subseteq \{v^{n}\mid n>0\}. \tag{$\ast$}
\end{equation*} 
Then $\sigma(\mathbf{ P_{z_{1},z_{2}}})=P_{z_{1},z_{2}}$ for all $z_{1}<z_{2}$ in $I$ and $\sigma(\mathbf{ M_{z_{1},z_{2}}^{s}})=M_{z_{1},z_{2}}^{s}$ for any $z_{1},z_{2}\in I$ such that $sz_{1}<z_{1}<z_{2}<sz_{2}$.
\end{Prop}
\begin{proof}
We have $\overline{\sigma(p)}=\sigma(\overline{p})$ for all $p\in \mathbf{ Z}[\Gamma]$. Moreover, the $\mathbf{ r}$-polynomials do not depend on the order, therefore we have $\sigma(\mathbf{ r}_{z_{1},z_{2}})=r_{z_{1},z_{2}}$ for any $z_{1},z_{2}\in W$.

We prove by induction on $\ell(z_{2})-\ell(z_{1})$ that $\sigma(\mathbf{ P}_{z_{1},z_{2}})=P_{z_{1},z_{2}}$ for all $z_{1}\leq z_{2}$ in I. If $\ell(z_{2})-\ell(z_{1})=0$ it is clear.

Assume that $\ell(z_{2})-\ell(z_{1})>0$. Applying $\sigma$ to \ref{P dependence} using the induction hypothesis yields 
\begin{align*}
\sigma(\overline{\mathbf{ P}}_{z_{1},z_{2}})-\sigma(\mathbf{ P}_{z_{1},z_{2}})&=\underset{z_{1}<z\leq z_{2}}{\sum}\sigma(\mathbf{r}_{z_{1},z})\sigma(\mathbf{ P}_{z,z_{2}})\\
&=\underset{z_{1}<z\leq z_{2}}{\sum}r_{z_{1},z}P_{z,z_{2}}.\\
\end{align*}
This relation and condition $(\ast)$ implies that $\sigma(\mathbf{ P}_{z_{1},z_{2}})=P_{z_{1},z_{2}}$.

Let $z_{1},z_{2}\in I$ and $s\in S$ such that $sz_{1}<z_{1}<z_{2}<sz_{2}$. We prove by induction on $\ell(z_{2})-\ell(z_{1})$ that $\sigma(\mathbf{ M}_{z_{1},z_{2}}^{s})=M_{z_{1},z_{2}}^{s}$.

Since $\mathbf{ M}_{z_{1},z_{2}}^{s}=\overline{\mathbf{ M}_{z_{1},z_{2}}^{s}}$, we have $\sigma(\mathbf{ M}_{z_{1},z_{2}}^{s})=\overline{\sigma(\mathbf{ M}_{z_{1},z_{2}}^{s})}$. Applying $\sigma$ to \ref{M dependence} (using $(\ast)$) we get 
$$ \sigma(\mathbf{ M}_{z_{1},z_{2}}^{s})+\underset{z;z_{1}< z<z_{2};sz<z}{\sum}P_{z_{1},z}M_{z,z_{2}}^{s}-v^{L(s)} P_{z_{1},z_{2}}
\in v^{-1}\mathbb{Z}[v^{-1}].$$
This relation implies that $ \sigma(\mathbf{ M}_{z_{1},z_{2}}^{s})=M_{z_{1},z_{2}}^{s}$.
Moreover we can see that if $\mathbf{ M}_{z_{1},z_{2}}^{s}\neq 0$ then $M_{z_{1},z_{2}}^{s}$ is a combination of pairwise different powers of $v$. Thus $M_{z_{1},z_{2}}^{s}\neq 0$.
\end{proof}

If condition $(\ast)$ is satisfied for all $s\in S$, then we can conclude that $x,z\in I$ satisfy $x\leftarrow_{L} z$ with respect to the total order $\Gamma_{+}$ if and only if they satisfy $x\leftarrow_{L} z$ with respect to the weight function $L$.

\section{The translations in an affine Weyl group}
We look more closely at a special set of elements in $W$, namely the translations. We keep the same notations as in Section 2.
\begin{Def}
\label{deftrans}
Let $u\in W$. We say that $u$ is a translation if there exists a vector $\vec{u}\neq 0$ such that $t_{\vec{u}}$, the translation by the vector $\vec{u}$, is in $\Omega$ and
$$uA_{0}=A_{0}t_{\vec{u}}.$$
Note that $t_{\vec{u}}$ is uniquely determined by $u$.
\end{Def}
Let $u\in W$ be a translation. Let  $B\in X$ and $\sigma\in\Omega$ be such that $B=A_{0}\sigma$. We have
$$uB=u(A_{0}\sigma)=A_{0}t_{\vec{u}}\sigma=A_{0}\sigma t_{\sigma(\vec{u})}=Bt_{\sigma (\vec{u})}.$$
Therefore $uB$ is a translate of $B$.

Recall that $h_{A}(z)$ for $z\in W$ and $A\in X$ is defined in \ref{hdef}.

The translations have the following properties.
\begin{Lem}
\label{basic}
Let $u\in W$ be a translation associated to $t_{\vec{u}}\in\Omega$.
$\ $\\
\begin{enumerate}
\item[(a)] Let $r_{1}\leq r_{2}\in \mathbb{N}^{*}$. We have
$$h_{A_{0}}(u^{r_{2}})\subset h_{A_{0}}(u^{r_{1}})\ \text{ and }\  h_{A_{0}}(u^{r_{2}})=t_{(r_{2}-r_{1})\vec{u}}(h_{A_{0}}(u^{r_{1}})).$$
\item[(b)] Let $r\in \mathbb{N}^{*}$. We have
$$z.u\Leftrightarrow z.u^{r}.$$
\end{enumerate}
\end{Lem}
\begin{proof} 
(a) Let $\alpha\in\Phi^{+}$ and $k_{\alpha}=\sg \vec{u},\check{\alpha}\sd$. Since $t_{\vec{u}}\in\Omega$, one can see that $k_{\alpha}\in\mathbb{Z}$. For any $r\in\mathbb{N}$, we have
$$rk_{\alpha}<\sg x,\check{\alpha}\sd<rk_{\alpha}+1\text{ for all $x\in u^{r}A_{0}$.}$$
Note that, if $k_{\alpha}=0$, there is no hyperplane of the form $H_{\alpha,m}$ ($m\in\mathbb{Z}$) which separates $A_{0}$ and $u^{r}A_{0}$.\\
Let $\varphi$ (resp. $\varphi^{+}$, $\varphi^{-}$) be the subset of $\Phi^{+}$ which consists of all positive roots $\beta$ such that $k_{\beta}\neq 0$ (resp. $k_{\beta}>0$,  $k_{\beta}<0$). For $\beta\in\varphi$, we define
$$H^{\beta}=
\begin{cases}
H_{\beta, rk_{\beta}} &\mbox{\text{ if $\beta\in\varphi^{+}$}},\\
H_{\beta, rk_{\beta}+1} &\mbox{\text{ if $\beta\in\varphi^{-}$.}}
\end{cases}
$$
Then, one can check that
\begin{equation*}
h_{A_{0}}(u^{r})=\bigcap_{\beta\in\varphi} E_{H^{\beta}}(u^{r}A_{0}). \tag{$\ast$}
\end{equation*}
Let $r_{1}\leq r_{2}\in \mathbb{N}^{*}$ and $\beta\in\varphi$. We suppose that $\beta\in\varphi^{+}$ (the case $\beta\in\varphi^{-}$ is similar).  We have
$$E_{H^{\beta}}(u^{r_{1}})=\{x\in V\ |\ \sg x,\check{\beta}\sd>r_{1}k_{\beta}\}$$
and
$$E_{H^{\beta}}(u^{r_{2}})=\{x\in V\ |\ \sg x,\check{\beta}\sd>r_{2}k_{\beta}\}.$$
Thus
$$E_{H^{\beta}}(u^{r_{2}})\subset  E_{H^{\beta}}(u^{r_{1}}) \quad\text{and}\quad      E_{H^{\beta}}(u^{r_{2}})=t_{(r_{2}-r_{1})\vec{u}}E_{H^{\beta}}(u^{r_{1}})$$
and the result follows using relation ($\ast$).

(b) The statement follows from (a) and Lemma \ref{hzone}.
\end{proof}


From now on and until the end of this section, we fix a translation $u\in W$ associated to $t_{\vec{u}}\in\Omega$. The orbit of $\vec{u}$ under $\Omega$ is finite. Indeed the group of linear transformations associated to the group of affine transformations $\Omega$ is finite, it is isomorphic to $\Omega_{0}$ which, in turn, is isomorphic to the Weyl group $W_{0}$ associated to the root system $\Phi$. Let
$$Orb_{\Omega}(\vec{u})=\{ \vec{u}_{1}=\vec{u},...,\vec{u}_{n}\}.$$
For $i\in\{1,...,n\}=[1,n]$, let $u_{i}\in W$ be such that $u_{i}A_{0}=A_{0}t_{\vec{u_{i}}}$, $v_{i}$ be the special point $t_{\vec{u}_{i}}(0)$ and $A_{0}t_{\vec{u}_{i}}=A_{v_{i}}$.

\begin{Lem}
\label{tool Th} 
\begin{enumerate}
\item[(a)] For any $i,j\in[1,n]$ we have $\ell(u_{i})=\ell(u_{j})$.
\item [(b)]Let $z_{1},z_{2}\in W$, $r\in\mathbb{N}^{*}$ and $i\in[1,n]$ be such that $z_{1}.u_{i}^{r}.z_{2}$. There exists $k,m\in[1,n]$ such that
$$z_{1}.u_{i}^{r}.z_{2}=z_{1}.z_{2}.u_{m}^{r}=u_{k}^{r}.z_{1}.z_{2}.$$
\item[(c)] Let $z_{1},z_{2}\in W$, $r\geq 1$ and $i\in [1,n]$. We have the following equivalence
$$z_{1}.u_{i}^{r}.z_{2} \Leftrightarrow z_{1}.u_{i}^{r+1}.z_{2}.$$
\end{enumerate}
\end{Lem}
\begin{proof}
(a) Let $A\in X$ and $A'$ be a translate of $A$ (by a translation in $\Omega$). Then the number of hyperplanes which separate $A$ and $A'$ is equal to the number of hyperplanes which separate $zA$ and $zA'$ for any $z\in W$.

Let $i,j\in[1,n]$, $\sigma\in \Omega$ and $z\in W$ be such that $\vec{u}_{i}\sigma=\vec{u}_{j}$ and $zA_{0}=A_{0}\sigma$. \\
We have
\begin{align*}
\ell(u_{i})&= | \{H\ |\ H \text{ separates } A_{0} \text{ and } A_{v_{i}}\}|\\
&=| \{H\ |\ H \text{ separates } z^{-1}A_{0} \text{ and } z^{-1}A_{v_{i}}\}|\\
&=| \{H\ |\  H \text{ separates } z^{-1}A_{0}\sigma \text{ and } z^{-1}A_{v_{i}}\sigma \}|.\\
\end{align*}
Since
$$z^{-1}A_{0}\sigma=A_{0}$$
and
$$z^{-1}A_{v_{i}}\sigma=z^{-1}A_{0}t_{\vec{u}_{i}}\sigma=z^{-1}A_{0}\sigma t_{\vec{u}_{j}}=
A_{v_{j}},$$
we obtain
\begin{align*}
\ell(u_{i})&=| \{H\ |\  H \text{ separates } z^{-1}A_{0}\sigma \text{ and } z^{-1}A_{v_{i}}\sigma \}|\\
&=| \{H\ |\  H \text{ separates } A_{0} \text{ and } A_{v_{j}} \}|\\
&=\ell(u_{j})
\end{align*}
as desired.\\ 
(b) Let $\sigma_{z_{1}},\sigma_{z_{2}}\in\Omega$ and $k,m\in[1,n]$ be such that 
$$z_{1}A_{0}=A_{0}\sigma_{z_{1}}\quad\ ,\quad \sigma^{-1}_{z_{1}}(\vec{u}_{i})=\vec{u}_{k}.$$
and
$$z_{2}A_{0}=A_{0}\sigma_{z_{2}}\quad , \quad \sigma_{z_{2}}(\vec{u}_{i})=\vec{u}_{m}.$$
We have
$$z_{1}.u_{i}^{r} .z_{2}A_{0}=A_{0}\sigma_{z_{1}}t_{r\vec{u}_{i}}\sigma_{z_{2}}=A_{0}t_{r\sigma_{z_{1}}^{-1}(\vec{u}_{i})}\sigma_{z_{1}}\sigma_{z_{2}}=u_{k}^{r}z_{1}z_{2}A_{0},$$
which implies that $z_{1}.u_{i}^{r} .z_{2}=u_{k}^{r}z_{1}z_{2}$. Now, since $\ell(u_{i})=\ell(u_{k})$, we must have $u_{k}^{r}.z_{1}.z_{2}$.\\
Similarly, one can show that $z_{1}.u_{i}^{r} .z_{2}=z_{1}.z_{2}.u_{m}^{r}$.

(c) The statement follows from (b) and Lemma \ref{basic} (b).

\end{proof}
We now state the main result of this section. 
\begin{Th}
\label{Th base}
Let $i,j\in[1,n]$ and $r_{1},r_{2}\in\bN^{*}$ be such that $i\neq j$. We have
$$h_{A_{0}}(u_{i}^{r_{1}})\cap h_{A_{0}}(u_{j}^{r_{2}})=\emptyset.$$
\end{Th}
\begin{proof}
According to Lemma \ref{basic} (a), to prove the theorem, it is enough to show that, for any $i\neq j\in[1,n]$, we have
$$h_{A_{0}}(u_{i})\cap h_{A_{0}}(u_{j})=\emptyset.$$
 Let 
$$\cF_{0}:=\{H\in\cF\ |\ 0\in H,\ v_{i}\notin H \text{ for all $i\in[1,n]$}\}.$$
Consider the connected component of 
$$V-\bigcup_{H\in\cF_{0}} H.$$
Since there exists $\sigma\in\Omega_{0}$ such that $\sigma(v_{i})=v_{j}$, there is a hyperplane which separate $v_{i}$ and $v_{j}$ and which contains $0$. Therefore $A_{v_{i}}$ and $A_{v_{j}}$ do not lie in the same connected component. For $i\in[1,n]$, let $\mathcal{C}_{i}$ be the connected component which contains $A_{v_{i}}$. To prove the theorem, it is enough to show that $h_{A_{0}}(u_{i})\subset\mathcal{C}_{i}$ for all $i\in[1,n]$.

Let $H$ be a wall of $\mathcal{C}_{i}$ and let $E_{H}(\mathcal{C}_{i})$ be the half-space defined by $H$ which contains $\mathcal{C}_{i}$. Since $0\in H$ and $v_{i}\notin H$, one can see that $H'=t_{\vec{u_{i}}}(H)\neq H$.  Thus either $H$ separates $A_{0}$ and $A_{v_{i}}$ or $H'$ does.

If $H$ separates $A_{0}$ and $A_{v_{i}}$ then, as $A_{v_{i}}\subset\mathcal{C}_{i}$, we must have $h_{A_{0}}(u_{i})\subset E_{H}(\mathcal{C}_{i})$.

Now, assume that $H$ does not separate $A_{0}$ and $A_{v_{i}}$. Let $\beta\in\Phi^{+}$ and $m\in\bZ$ be such that $H=H_{\beta,0}$ and $H'=H_{\beta,m}$. In that case we have
$$ A_{0},A_{v_{i}}\in E_{H}(\mathcal{C}_{i})=\{x\in V\mid \mathcal{h}x,\check{\beta}\mathcal{i}>0\}.$$
Thus one can see that we must have $m>0$. Let $E_{H'}(A_{v_{i}})$ be the half-space defined by $H'$ which contains $A_{v_{i}}$. We have
$$E_{H'}(A_{v_{i}})=\{x\in V\ |\ \sg x,\check{\beta}\sd>m\}$$
and 
$$h_{A_{0}}(u_{i})\subset E_{H'}(A_{v_{i}})\subset E_{H}(\mathcal{C}_{i}).$$
We have shown that for every wall of $\mathcal{C}_{i}$, $h_{A_{0}}(u_{i})$ lies on the same side of this wall as $\mathcal{C}_{i}$, thus $h_{A_{0}}(u_{i})\subset \mathcal{C}_{i}$ as required.
\end{proof}

\begin{Cor} 
(of Theorem \ref{Th base})\\
\label{Th utile}
\begin{enumerate}
\item[(a)] Let $z,z'\in W$, $r\in\mathbb{N}^{*}$, $m\in\mathbb{N}$ and $i,j\in[1,n]$. We have
$$z.u_{i}^{r}=z'.u_{j}^{r+m}\Longrightarrow i=j \text{ and }z=z'.u_{j}^{m}.$$
\item[(b)] Let $z_{1},z_{2},z'_{1},z'_{2}\in W$, $r\in\mathbb{N}^{*}$, $m\in\mathbb{N}$ and $i,j\in[1,n]$. For all $k\geq 0$ we have
$$z_{1}.u_{i}^{r}.z_{2}=z'_{1}.u_{j}^{r+m}.z'_{2} \Leftrightarrow z_{1}.u_{i}^{r+k}.z_{2}=z'_{1}.u_{j}^{r+k+m}.z'_{2}.$$
\end{enumerate}
\end{Cor}
\begin{proof}
(a) We have $z.u_{i}^{r}A_{0}\in h_{A_{0}}(u_{i}^{r})$ and $z'.u_{j}^{r+m}=z'.u_{j}^{m}.u_{j}^{r}\in h_{A_{0}}(u_{j}^{r})$. Since $z.u_{i}^{r}=z'.u_{j}^{r+m}$, applying Theorem \ref{Th base} yields $i=j$. The result follows.\\

(b) The statement follows from Lemma \ref{tool Th} and (a).
\end{proof}

\section{Isomorphism of intervals and equalities of Kazhdan-Lusztig polynomials}
Let $u\in W$ be a translation associated to $t_{\vec{u}}\in\Omega$ and let $M=\ell(u)$. One can easily see that $M\geq 2$.
Let $$Orb_{\Omega}(\vec{u})=\{\vec{u}_{1}=\vec{u},...,\vec{u}_{n}\}.$$
Finally, for $i\in[1,n]$, let $u_{i}\in W$ be such that $u_{i}A_{0}=A_{0}t_{\vec{u}_{i}}$. In this section we want to prove
\begin{thm1}
Let $z,z'\in W$, $r\in \mathbb{N}^{*}$ and $i,j\in[1,n]$ be such that $z.u_{i}^{r}$ and $z'.u_{j}^{r}$. Let $N=\ell(z')-\ell(z)$. Then for $r>N(M+1)$ and for any $k\geq 0$ we have
$$\mathbf{ P}_{z.u_{i}^{r},z'.u_{j}^{r}}=\mathbf{ P}_{z.u_{i}^{r+k},z'.u_{j}^{r+k}}$$
and if there exists $s\in S$ which satisfies $sz.u_{i}^{r}<z.u_{i}^{r}<z'.u_{j}^{r}<s.z'.u_{j}^{r}$, then
$$\mathbf{ M}_{z.u_{i}^{r},z'.u_{j}^{r}}^{s}=\mathbf{ M}_{z.u_{i}^{r+k},z'.u_{j}^{r+k}}^{s}\ .$$
\end{thm1} 
Our first task is to construct an isomorphism from the Bruhat interval $[z.u_{i}^{r},z'.u_{j}^{r}]$ to $[z.u_{i}^{r+k},z'.u_{j}^{r+k}]$ and then to show that the corresponding $\mathbf{ r}$-polynomials are equal.\\
\begin{Lem}
\label{proof prop}
Let $z,y\in W$, $r\in \mathbb{N}^{*}$ and $i\in [1,n]$ be such that $z.u_{i}^{r}$. Let $N=\ell(z.u_{i}^{r})-\ell(y)$. Then for $r>N$ we have
\begin{eqnarray*}
y\leq  z.u_{i}^{r}&\Leftrightarrow& \exists z_{1},z_{2}\in W, n_{1},n_{2}\in \mathbb{N} \text{ such that } z_{1}.u_{i}^{r-N}.z_{2}=y \\
& & \ \ \ \ \ z_{1}\leq z.u_{i}^{n_{1}},z_{2}\leq u_{i}^{n_{2}} \text{ and $n_{1}+n_{2}=N$.}  
\end{eqnarray*}
Furthermore, there exists a unique $z_{y}\in W$ and $m\in [1,n]$ such that $y=z_{y}.u_{m}^{r-N}$.
\end{Lem}
\begin{proof}
$``\Leftarrow''$ is clear.\\
$``\Rightarrow''$ We proceed by induction on $N$. \\
If $N=0$, it's clear.\\
Let $N>0$. There exists $y'\in W$ such that $y\leq y'\leq z.u_{i}^{r}$ and
$$\ell(z.u_{i}^{r})-\ell(y')=N-1\ \text{ and }\ \ell(y')-\ell(y)=1 .$$
Applying the inductive assumption yields
\begin{eqnarray*}
& & \exists z'_{1}, z'_{2}\in W, n'_{1}, n'_{2}\in \mathbb{N} \text{ such that } z'_{1}.u_{i}^{r-N+1}.z'_{2}=y'\\
& &\ \ \ z'_{1}\leq z.u_{i}^{n'_{1}}, z'_{2}\leq u_{i}^{n'_{2}}, n'_{1}+n'_{2}=N-1\ .\\
\end{eqnarray*}
Let
$$y'=s_{p}...s_{m+1}(s_{m}...s_{k})^{r-N+1}s_{k-1}...s_{1}\qquad (p\geq m\geq k\geq 1)$$
be a reduced expression of $y'$ such that
$$z'_{1}=s_{p}...s_{m+1},\ u_{i}=s_{m}...s_{k}\text{ and }z'_{2}=s_{k-1}...s_{1}\ .$$
We know that $y$ can be obtained by dropping a simple reflection $s\in S$ in a reduced expression of $y'$. If there exists $l\in\mathbb{N}$ such that 
$$y=s_{p}...\hat{s_{l}}...s_{k}(s_{m}...s_{k})^{r-N}s_{k-1}...s_{1}\qquad (p\geq l\geq k)$$
or
$$y=s_{p}...s_{m+1}(s_{m}...s_{k})^{r-N}s_{m}...\hat{s_{l}}...s_{1}\qquad (m\geq l\geq 1)$$
(where $\hat{s}$ means that we have dropped $s$) the result is straightforward.\\
Now assume that there exists $l_{1},l_{2}\in\mathbb{N}^{*}$ such that $l_{1}+l_{2}=r-N$ and
$$y=z'_{1}.u_{i}^{l_{1}}.\hat{u}_{i}.u_{i}^{l_{2}}.z'_{2}\ .$$
where $\hat{u}_{i}$ is obtained by dropping a simple reflection in $s_{m}...s_{k}$.\\
Let $j\in[1,n]$ such that $u_{i}^{l_{1}}.\hat{u}_{i}=\hat{u}_{i}.u_{j}^{l_{1}}$. We have $y=z'_{1}.\hat{u}_{i}.u_{j}^{l_{2}}.u_{i}^{l_{1}}.z'_{2}$ which implies that $u_{j}^{l_{2}}.u_{i}^{l_{1}}$. Furthermore, we have
$$u_{j}^{l_{2}}.u_{i}^{l_{1}}A_{0}=A_{0}t_{l_{2}\vec{u}_{j}}t_{l_{1}\vec{u}_{i}}=A_{0}t_{l_{1}\vec{u}_{i}}t_{l_{2}\vec{u}_{j}}=u_{i}^{l_{1}}.u_{j}^{l_{2}}A_{0}.$$
Applying Corollary \ref{Th utile}, we get $i=j$. Thus
$$y=z'_{1}.u_{i}^{l_{1}}.\hat{u}_{i}.u_{i}^{l_{2}}.z'_{2}=z'_{1}.\hat{u}_{i}.u_{i}^{r-N}.z'_{2}.$$
Let
\begin{eqnarray*}
&z_{1}=z'_{1}.\hat{u}_{i} & n_{1}=n'_{1}+1\\
&z_{2}=z'_{2} & n_{2}=n'_{2}.
\end{eqnarray*}
Then one can check that $z_{1}=z'_{1}.\hat{u}_{i}\leq z.u_{i}^{n_{1}}$ and $z_{2}\leq u_{i}^{n_{2}}$. Thus we get the result by induction.

Let $m\in [1,n]$ be such that $y=z_{1}.u_{i}^{r-N}.z_{2}=z_{1}.z_{2}.u_{m}^{r-N}$. Let $z_{y}=z_{1}.z_{2}$. Assume that there exists $w\in W$ and $k\in [1,n]$ such that $y=w.u_{k}^{r-N}$. By Corollary \ref{Th utile}, we have $k=m$ and $w=z_{y}$, which concludes the proof.
\end{proof}

\begin{Lem}
\label{Th isom}
Let $z,z'\in W$, $i,j\in [1,n]$ and $r_{1},r_{2}\in\mathbb{N}^{*}$ be such that $z.u_{i}^{r_{1}}$ and $z'.u_{j}^{r_{2}}$. Let $N=\ell(z'.u_{j}^{r_{2}})-\ell(z.u_{i}^{r_{1}})$. Then for $r_{2}>N$ and for any $k\geq 0$ we have
$$z.u_{i}^{r_{1}} \leq z'.u_{j}^{r_{2}} \Leftrightarrow z.u_{i}^{r_{1}+k}\leq z'.u_{j}^{r_{2}+k}.$$
\end{Lem}
\begin{proof}
Applying the previous lemma and Corollary \ref{Th utile} yields the following equivalences, for any $k\geq 0$
\begin{align*}
&z.u_{i}^{r_{1}}\leq z'.u_{j}^{r_{2}} \\
&\Leftrightarrow \exists z_{1},z_{2}\in W, n_{1},n_{2}\in \mathbb{N} \text{ such that } z_{1}.u_{j}^{r_{2}-N}.z_{2}=z.u_{i}^{r_{1}}\\
&  \ \ \ \ \ z_{1}\leq z'.u_{j}^{n_{1}},z_{2}\leq u_{j}^{n_{2}}  \text{ and $n_{1}+n_{2}=N$}  \\
&\Leftrightarrow  \exists z_{1},z_{2}\in W, n_{1},n_{2}\in \mathbb{N} \text{ such that } z_{1}.u_{j}^{r_{2}-N+k}.z_{2}=z.u_{i}^{r_{1}+k}, \\
& \ \ \ \ \ \  z_{1}\leq z'.u_{j}^{n_{1}},z_{2}\leq u_{j}^{n_{2}}  \text{ and $n_{1}+n_{2}=N$}\\
&\Leftrightarrow zu_{i}^{r_{1}+k}\leq z'u_{j}^{r_{2}+k}.
\end{align*}
\end{proof}
\begin{Prop}
\label{isomorphism}
Let $z,z'\in W$,  $i,j\in [1,n]$ and $r\in\mathbb{N}^{*}$ be such that $z.u_{i}^{r}$ and $z'.u_{j}^{r}$. Let $N=\ell(z'.u_{j}^{r})-\ell(z.u_{i}^{r})$. Then for $r>N$ and for any $k\geq 0$,  the Bruhat interval
$$I_{1}=[z.u_{i}^{r},z'.u_{j}^{r}]=\{y\in W\mid z.u_{i}^{r}\leq y\leq z'.u_{j}^{r}\}$$
is isomorphic to $I_{2}=[z.u_{i}^{r+k},z'.u_{j}^{r+k}]$.
\end{Prop}
\begin{proof}
Let $y\in I_{1}$ and $N_{y}=\ell(z'.u_{j}^{r})-\ell(y)$. There exists a unique $z_{y}\in W$ and $m\in\mathbb{N}$ such that $y=z_{y}.u_{m}^{r-N_{y}}$.\\
Let
$$
\begin{array}{cccc}
\varphi: &I_{1}&\longrightarrow &I_{2}\\
&z_{y}.u_{m}^{r-N_{y}} &\longmapsto &z_{y}.u_{m}^{r+k-N_{y}}.
\end{array}
$$
We need to show that $\varphi$ is an isomorphism of Bruhat interval.\\
Let $y'\leq y\in I_{1}$. Let $N_{y'}=\ell(z'.u_{j}^{r})-\ell(y')$. There exists a unique $z_{y'}\in W$ and $m'\in\mathbb{N}$ such that $y=z_{y'}.u_{m'}^{r-N_{y'}}$. One can check that we can apply Lemma \ref{Th isom}, we obtain
\begin{align*}
&\quad z.u_{i}^{r}\leq y=z_{y}.u_{m}^{r-N_{y}}\leq y'=z_{y'}.u_{m'}^{r-N'_{y}} \leq z'.u_{j}^{r}\\
\Longleftrightarrow&\quad z.u_{i}^{r+k}\leq \varphi(y)=z_{y}.u_{m}^{r-N_{y}+k}\leq \varphi(y')=z_{y'}.u_{m'}^{r-N_{y'}+k} \leq z'.u_{j}^{r+k}.
\end{align*}
By Corollary \ref{Th utile} we see that $\varphi$ is injective. One can easly check that $\varphi$ is surjective. The result follows.
\end{proof}

The next step is to show that the corresponding $\mathbf{r}$-polynomials are equal. Let $\Gamma$ be an abelian group together with a total order specified by $\Gamma_{+}$. Let $\{v_{s}\mid s\in S\}\subset \Gamma_{+}$ be the set of parameters and $\xi_{s}=v_{s}-v_{s}^{-1}$ (see Section 3 for details). 

Let $y,w\in W$ and $s\in S$ such that $sw<w$. Recall that the $\mathbf{ r}$-polynomials satisfy $\mathbf{ r}_{y,w}=0$ unless $y\leq w$, $\mathbf{ r}_{y,y}=1$ and the  recursive relation
\begin{equation*}
\label{rec}
\mathbf{ r}_{y,w}=
\begin{cases}
\mathbf{ r}_{sy,sw},&\text{if}\ sy<y,\\
\mathbf{ r}_{sy,sw}+(v_{s}-v_{s}^{-1})\mathbf{ r}_{y,sw},&\text{if}\ sy>y\ .\\
\end{cases}
\end{equation*}
\begin{Prop}
\label{rpol}
Let $z,z'\in W$, $i,j\in [1,n]$ and $r\in \mathbb{N}^{*}$ be such that $z.u_{i}^{r}$ and $z'.u_{j}^{r}$. Let $N=\ell(z')-\ell(z)$. Then for $r>NM$ and for any $k\geq 0$ we have
$$\mathbf{r}_{z.u_{i}^{r},z'.u_{j}^{r}}=\mathbf{r}_{z.u_{i}^{r+k},z'.u_{j}^{r+k}}.$$ 
\end{Prop}
\begin{proof}
We proceed by induction on $N$.\\
If $i=j$ or if $N<0$ then the result is obvious.\\
If $N=0$, since
$$\mathbf{r}_{z.u_{i}^{r},z'.u_{j}^{r}}=\delta_{z.u_{i}^{r},z'.u_{j}^{r}}\ \text{ and } \mathbf{r}_{z.u_{i}^{r+k},z.u_{j}^{r+k}}=\delta_{z.u_{i}^{r+k},z.u_{j}^{r+k}}$$ 
the result follows from Corollary \ref{Th utile}.\\
Let $N\geq 1$ and $i\neq j$. Note that in this case $r>M$.\\
Let $u_{j}=s_{M}...s_{1}$ be a reduced expression. There exists $1\leq l\leq M$ such that $(z.u_{i}^{r}s_{1}...s_{l-1})s_{l}>z.u_{i}^{r}s_{1}...s_{l-1}$. Indeed, if not, then $z.u_{i}^{r}=y.u_{j}$ for some $y\in W$. By Corollary \ref{Th utile}, this implies that $i=j$, but we assumed that $i\neq j$. Let $l\in[1...m]$ be the smallest element with this property. The minimality of $l$ implies that $\ell(z.u_{i}^{r}s_{1}...s_{l-1})=\ell(z.u_{i}^{r})-(l-1)$.\\
One can see that $zu_{i}^{r}s_{1}...s_{l-1}\leq z.u_{i}^{r}$. Let $y,w\in W$ and $m,q,p\in[1,n]$ be such that
\begin{align*}
zu_{i}^{r}s_{1}...s_{l-1}&=y.u_{m}^{r-l+1}\\
y.u_{m}^{r-l+1}.s_{l}&=y.s_{l}.u_{p}^{r-l+1}\\
z'u_{j}^{r}s_{1}...s_{l}&=z'.u_{j}^{r-1}.s_{M}...s_{l+1}=w.u_{q}^{r-1}\\
\end{align*}
By Corollary \ref{Th utile} we see that
\begin{align*}
zu_{i}^{r+k}s_{1}...s_{l-1}&=y.u_{m}^{r+k-l+1}\\
y.u_{m}^{r+k-l+1}.s_{l}&=y.s_{l}.u_{p}^{r+k-l+1}\\
z'u_{j}^{r+k}s_{1}...s_{l}&=z'.u_{j}^{r+k-1}.s_{M}...s_{l+1}=w.u_{q}^{r+k-1}.
\end{align*}
Applying the recursive formula for the $\mathbf{r}$-polynomials, we obtain
\begin{align*}
\mathbf{r}_{z.u_{i}^{r},z.u_{j}^{r}}&=\mathbf{r}_{zu_{i}^{r}s_{1}...s_{l-1},z_{2}u_{j}^{r}s_{1}...s_{l-1}}\\
&=\mathbf{r}_{y.u_{m}^{r-l+1},z_{2}u_{j}^{r}s_{1}...s_{l-1}}\\
&=\mathbf{r}_{ys_{l}.u_{p}^{r-l+1},w.u_{q}^{r-1}}+\xi_{s_{l}}\mathbf{r}_{y.u_{m}^{r-l+1},w.u_{q}^{r-1}}\\
\end{align*}
and
\begin{align*}
\mathbf{r}_{z.u_{i}^{r+k},z.u_{j}^{r+k}}&=\mathbf{r}_{zu_{i}^{r+k}s_{1}...s_{l-1},z_{2}u_{j}^{r+k}s_{1}...s_{l-1}}\\
&=\mathbf{r}_{y.u_{m}^{r-l+1+k},z_{2}u_{j}^{r+k}s_{1}...s_{l-1}}\\
&=\mathbf{r}_{ys_{l}.u_{p}^{r-l+1+k},w.u_{q}^{r-1+k}}+\xi_{s_{l}}\mathbf{r}_{y.u_{m}^{r-l+1+k},w.u_{q}^{r-1+k}}\ \ .\\
\end{align*}
Therefore to prove the theorem it is enough to show that
\begin{align*}
\mathbf{r}_{ys_{l}.u_{p}^{r-l+1},w.u_{q}^{r-1}}&=\mathbf{r}_{ys_{l}.u_{p}^{r-l+1+k},w.u_{q}^{r-1+k}},\\ 
\mathbf{r}_{y.u_{m}^{r-l+1},w.u_{q}^{r-1}}&=\mathbf{r}_{y.u_{m}^{r-l+1+k},w.u_{q}^{r-1+k}}.
\end{align*}
If $l=1$ we have
\begin{align*}
\mathbf{r}_{ys_{l}.u_{p}^{r},w.u_{q}^{r-1}}&=\mathbf{r}_{(ys_{l}.u_{p}).u_{p}^{r-1},w.u_{q}^{r-1}},\\
\mathbf{r}_{y.u_{m}^{r},w.u_{q}^{r-1}}&=\mathbf{r}_{(y.u_{m}).u_{m}^{r-1},w.u_{q}^{r-1}}
\end{align*}
and
\begin{align*}
\ell(w)-\ell(ys_{l}u_{p})&=N-2,\\
\ell(w)-\ell(yu_{m})&=N-1.
\end{align*}
Moreover $r-1>0$ (we have seen that $r>M\geq 2$) and
$$r-1>MN-1\geq MN-N=M(N-1)>M(N-2).$$
Therefore in both cases we can apply the induction hypothesis which yields the desired equalities.\\
If $l>1$, we have
\begin{align*}
\mathbf{r}_{ys_{l}u_{p}^{r-l+1},w.u_{q}^{r-1}}&=\mathbf{r}_{ys_{l}u_{p}^{r-l+1},w.u_{q}^{l-2}.u_{q}^{r-l+1}},\\
\mathbf{r}_{y.u_{m}^{r-l+1},w.u_{q}^{r-1}}&=\mathbf{r}_{y.u_{m}^{r-l+1},w.u_{q}^{l-2}u_{q}^{r-l+1}}
\end{align*}
and
\begin{align*}
\ell(w.u_{q}^{l-2})-\ell(ys_{l})&=N-2,\\
\ell(w.u_{q}^{l-2})-\ell(y)&=N-1.
\end{align*}
Moreover $r-l+1>0$ and
$$r-l+1>r-M>M(N-1)>M(N-2)$$
and once more the induction hypothesis gives the desired equalities.
\end{proof}

We are now ready to prove Theorem \ref{translation}.
\begin{proof}
The intervals $I_{1}=[z.u_{i}^{r},z'.u_{j}^{r}]$, $I_{2}=[z.u^{r+k},z'.u^{r+k}]$ are isomorphic with respect to the Bruhat order via $\varphi$ (as defined in Proposition \ref{isomorphism}). \\
Let $y_{1}, y_{2}\in I$ ($\ell(y_{1})\leq \ell(y_{2})$), $z_{1},z_{2}\in W$, $N_{1},N_{2}\in\bN$ and $m_{1},m_{2}\in[1,n]$ be such that 
$$N_{1}=\ell(z'.u^{r+k})-\ell(y_{1})\quad\text{and}\quad y_{1}=z_{1}u_{m_{1}}^{r-N_{1}},$$
$$N_{2}=\ell(z'.u^{r+k})-\ell(y_{2})\quad\text{and}\quad y_{2}=z_{2}u_{m_{2}}^{r-N_{2}}.$$
We have
$$r-N_{1}\geq r-N>N(M+1)-N=MN\geq M(\ell(y_{2})-\ell(y_{1})).$$
Thus, by Proposition \ref{rpol}, we obtain
\begin{align*}
r_{y_{1},y_{2}}&=r_{z_{1}u_{m_{1}}^{r-N_{1}},z_{2}u_{m_{2}}^{r-N_{2}}}\\
&=r_{z_{1}u_{m_{1}}^{r-N_{1}},z_{2}u_{m_{2}}^{N_{1}-N_{2}}u_{m_{2}}^{r-N_{1}}}\\
&= r_{z_{1}u_{m_{1}}^{r+k-N_{1}},z_{2}u_{m_{2}}^{N_{1}-N_{2}}u_{m_{2}}^{r+k-N_{1}}}\\
&=r_{\varphi(y_{1}),\varphi(y_{2})}.
\end{align*}
Therefore, by remark \ref{P,M dependence}, we get the result.
\end{proof}

\section{Application to $\tilde{G}_{2}$}
Our aim is to prove Theorem \ref{finite}.\\ 
Throughout this section, let $W$ be an affine Weyl group of type $\tilde{G}_{2}$, with presentation as follows
$$W:=\sg s_{1},s_{2},s_{3}\ |\ (s_{1}s_{2})^{6}=1,(s_{2}s_{3})^{3}=1,(s_{1}s_{3})^{2}=1\sd.$$

The generators $s_{2}$ and $s_{3}$ are conjugate in W, thus a weight function $L$ on $W$ is uniquely determined by
$$L(s_{1})=a \ \text{  and  } \ L(s_{2})=L(s_{3})=b\quad a,b\in\bN^{*}$$ 
We shall denote such a weight function by $L=L_{a,b}$.

Let $Q,q$ be independant indeterminates over $\mathbb{Z}$ and consider the abelian group
$$\Gamma=\{Q^{i}q^{j}\mid i,j\in\mathbb{Z}\}.$$
Let $v$ be another indeterminate. We have a ring homomorphism
$$\sigma_{a,b}: \mathbb{Z}[\Gamma]\rightarrow \mathbb{Z}[v,v^{-1}],\ \ \ \ Q^{i}q^{j} \rightarrow v^{ai+bj}.$$
We will need the following lemma.
\begin{Lem}
\label{calcul}
Let $y<w\in W$, $I=[y,w]$ and $s\in S$ such that $$sy<y<w<sw.$$
\begin{enumerate}
\item Consider the total order given by
$$\Gamma^{1}_{+}=\{Q^{i}q^{j}\mid i>0, j\in\mathbb{Z}\}\cup\{q^{i}\mid i>0\}.$$
Suppose that, for $c,d\in \mathbb{N}^{*}$, we have
$$\Gamma_{+}^{s}(I)\subseteq \{q^{j}\mid j>0\}\cup \{Q^{i}q^{j}\mid i>0, ci+dj\geq 0\}.$$
Then condition ($\ast$) in Proposition \ref{L} holds for any $\sigma_{a,b}$ such that $a/b>c/d$.\\
Furthermore, if $\mathbf{ M}_{y,w}^{s}\neq 0$, then for any weight functions $L_{a,b}$ such that $a/b>c/d$, we have $M_{y,w}^{s}\neq 0$.\\
\item Let $c\geq d\in\mathbb{N}^{*}$. Consider the total order given by
$$\Gamma^{2}_{+}=\{Q^{i}q^{j}\mid ci+dj>0\} \cup \{ Q^{dj}q^{-cj}\mid j>0\}.$$
Suppose that we have, for some $e>c/d\in \mathbb{Q}_{>0}$
$$\Gamma_{+}^{s}(I)\subseteq \{ q^{j}\mid j>0\}\cup \{Q^{i}q^{j}\mid i>0,i+j\geq 0\}$$
$$\cup \{Q^{i}q^{j}\mid j>-i>0, -j/i\geq e\}\cup \{Q^{i}q^{j}\mid -j>i>0,-j/i\leq c/d\}.$$
Then condition ($\ast$) in Proposition \ref{L} holds for any $\sigma_{a,b}$ such that $e>a/b>c/d$.\\
Furthermore, if $\mathbf{ M}_{y,w}^{s}\neq 0$, then for any weight functions $L_{a,b}$ such that $e>a/b>c/d$, we have $M_{y,w}^{s}\neq 0$.
\end{enumerate}
\end{Lem}
\begin{proof}
We prove 1. Let $i,j\in \mathbb{Z}$ such that $Q^{i}q^{j}\in \Gamma^{s}_{+}(I)$. We must show that $ai+bj>0$ provided that $a/b>c/d$.\\
 If $i=0$ then $j>0$ and $ai+bj=bj>0$. \\
 If $i>0$ and $ci+dj\geq 0$ then
 $$ai+bj=b(ia/b+j)>b(ic/d+j)\geq 0$$ as required.
 
We prove 2. Let $i,j\in \mathbb{Z}$ such that $Q^{i}q^{j}\in \Gamma^{s}_{+}(I)$. We must show that $ai+bj>0$ provided that $e>a/b>c/d$.\\
If $i=0$ then $j>0$ and $ai+bj=bj>0$. \\
If $i>0$ and $i+j\geq 0$ then $$ai+bj>b(c/d)i+bj=b(ic/d+j)>b(i+j)\geq 0.$$
If $j>-i>0$ and $-j/i\geq e$ then $$ai+bj=bj((a/b)(i/j)+1)>bj(ei/j+1)\geq 0.$$
Finally, if $-j>i>0$ and $-j/i\leq c/d$ then $$ai+bj=ai(1+(b/a)j/i)>ai(1+(d/c)(i/j))\geq 0$$ as required.\\
\end{proof}
Note that, in the situation of the above lemma, we will always have
$a>b$. But similar results also hold for $b>a$.\\

In order to prove Theorem \ref{finite}, we will proceed as follows.\\
Using Proposition \ref{L}, Theorem \ref{translation} and the previous lemma, we will find a collection of non-zero $M$-polynomials. We will then find some infinite sets such that each of these sets is included in a left cell for any choice of parameters and such that all the elements of $W$ lie in one of these sets except for a finite number. Then, we can conclude that the number of left cells is finite for any choice of parameters and that there is a finite number of distinct decompositions of $W$ into left cells.

We have developped some programs in GAP3 \cite{GAP} which, given an interval $I$, $s\in S$ and a monomial order on $\Gamma$, compute the following data
\begin{enumerate}
\item The Kazhdan-Lusztig polynomials $\mathbf{P}_{y,w}$ for all $y,w\in I$,
\item $\mathbf{M}_{y,w}^{s}$ for all $y,w\in I$ such that $sy<y<w<sw$,
\item $\underset{z;z_{1}\leq z<z_{2};sz<z}{\sum}\mathbf{ P}_{z_{1},z}\mathbf{ M}_{z,z_{2}}^{s}-v_{s}\mathbf{ P}_{z_{1},z_{2}}$ for all $z_{1},z_{2}\in I$,
\end{enumerate}
so that we can compute the set $\Gamma_{+}^{s}(I)$ as described in Proposition \ref{L}.\\

We work with the geometric representation of $W$ as described in Section 2. 
Let $u_{1}=s_{1}s_{2}s_{1}s_{2}s_{3}s_{1}s_{2}s_{1}s_{2}s_{3}\in W$. One can check that $u_{1}$ is a translation. Let
\begin{align*}
\Pi&=\{e,\ s_{3},\ s_{2}s_{3},\ s_{1}s_{2}s_{3},\ s_{2}s_{1}s_{2}s_{3},\\
& \ \ \ \ \ s_{3}s_{2}s_{1}s_{2}s_{3},\ s_{1}s_{2}s_{1}s_{2}s_{3},\ s_{3}s_{1}s_{2}s_{1}s_{2}s_{3},\ s_{2}s_{3}s_{1}s_{2}s_{1}s_{2}s_{3},\\
& \ \ \ \ \ s_{1}s_{2}s_{3}s_{1}s_{2}s_{1}s_{2}s_{3},\ s_{2}s_{1}s_{2}s_{3}s_{1}s_{2}s_{1}s_{2}s_{3},\ s_{3}s_{2}s_{1}s_{2}s_{3}s_{1}s_{2}s_{1}s_{2}s_{3}\}\\
W_{1}&=\{e,s_{1}, s_{1}s_{2},s_{1}s_{2}s_{1},s_{1}s_{2}s_{1}s_{2},s_{1}s_{2}s_{1}s_{2}s_{1}\}\\
&=\{w_{1},w_{2},w_{3},w_{4},w_{5},w_{6}\}
\end{align*}
and $y=s_{3}s_{2}s_{1}s_{2}s_{3}s_{1}s_{2}s_{1}s_{2}s_{3}$.\\
For $w_{i}\in W_{1}$ let $\sigma_{w_{i}}\in \Omega$ such that $w_{i}A_{0}=A_{0}\sigma_{w_{i}}$. One can check that
$$Orb_{\Omega}(\vec{u_{1}})=\{\vec{u}_{1}\sigma_{w_{1}},...,\vec{u}\sigma_{w_{6}}\}=\{\vec{u}_{1},...,\vec{u}_{6}\}.$$
For $1\leq i\leq 6$, let
$$I_{1,i}^{r}=[s_{3}.u^{r}.w_{i},y.u^{r}.w_{i}]=[s_{3}.w_{i}.u_{i}^{r},y.w_{i}.u_{i}^{r}]=[x_{1,i}^{r},y_{1,i}^{r}],$$
\begin{align*}
I_{2,i}^{r}&=[s_{3}s_{1}s_{2}s_{1}s_{2}s_{3}u^{r}.w_{i},y.s_{1}s_{2}s_{1}s_{2}s_{3}u^{r}.w_{i}]\\
&=[s_{3}s_{1}s_{2}s_{1}s_{2}s_{3}.w_{i}.u_{i}^{r},y.s_{1}s_{2}s_{1}s_{2}s_{3}.w_{i}.u_{i}^{r}]\\
&=[x_{2,i}^{r},y_{2,i}^{r}].
\end{align*}
Let $k=1,2$. According to the results in the previous section, we know that for $r$ large enough, the intervals $I_{k,i}^{r}$ and $I_{k,i}^{r+n}$ are isomorphic for any $n\geq 0$ and that the corresponding $\mathbf{r}$-polynomials are equal. However, using our GAP3 program, we see that this is true for $r\geq 6$. \\
We want to show that, for $r\geq 6$, $M_{x_{k,i}^{r},y_{k,i}^{r}}^{s_{1}}\neq 0$ for any choice of parameters. By the previous remark, it is enough to show that $M_{x_{k,i}^{6},y_{k,i}^{6}}^{s_{1}}\neq 0$.\\
We give some details for the computation of $M_{x_{1,1}^{6},y_{1,1}^{6}}^{s_{1}}$.\\
Let $x=x_{1,1}^{6}$, $y=y_{1,1}^{6}$ and $I=I_{1,1}^{6}$.\\
Consider the total order given by
$$\Gamma_{+}=\{Q^{i}q^{j}\mid i>0, j\in\mathbb{Z}\}\cup\{q^{i}\mid i>0\}.$$
Using our GAP3 program to calculate the set $\Gamma_{+}^{s_{1}}(I)$, we find $\mathbf{ M}_{x,y}^{s_{1}}=1$ and
$$\Gamma_{+}^{s}(I)\subseteq \{q^{j}\mid j>0\}\cup \{Q^{i}q^{j}\mid 3i+j\geq 0\}.$$
Therefore, applying Lemma \ref{calcul}, we see that, for any parameter $a,b$ such that $a/b>3$, we have 
$$M_{x,y}^{s_{1}}=\sigma_{a,b}(\mathbf{ M}_{x,y}^{s_{1}})\quad\text{and}\quad\mathbf{ M}_{x,y}^{s_{1}}\neq 0\Longrightarrow M_{x,y}^{s_{1}}\neq 0.$$
In order to deal with weight functions $L_{a,b}$ such that $a/b<3$, we proceed as follows. Let
$$\mathcal{E}=\{x\in\mathbb{Q}_{>0}\mid x=\pm j/i\text{ where } j<0,\ i\neq 0, Q^{i}q^{j}\in \Gamma_{+}^{s}(I) \}.$$
The largest element of $\mathcal{E}$ below $3$ is $2$. Thus we consider the weight functions $L_{a',b'}$ where $b'/a'>2$.\\
Consider the total order given by
$$\Gamma_{+}=\{Q^{i}q^{j}\mid 2i+j>0\} \cup \{ Q^{j}q^{-2j}\mid j>0\}.$$
Computing $\Gamma_{+}^{s}(I)$ gives $\mathbf{ M}_{x,y}^{s_{1}}=1$ and
$$\Gamma_{+}^{s}(I)\subseteq \{ q^{j}\mid j>0\}\cup \{Q^{i}q^{j}\mid i>0,i+j\geq 0\}$$
$$\cup \{Q^{i}q^{j}\mid j>-i>0, -j/i\geq 3\}\cup \{Q^{i}q^{j}\mid -j>i>0,-j/i\leq 2\}.$$
Therefore, for any parameter $a,b$ such that $3>a/b>2$, we have 
$$M_{x,y}^{s_{1}}=\sigma_{a,b}(\mathbf{ M}_{x,y}^{s_{1}})\quad\text{and}\quad\mathbf{ M}_{x,y}^{s_{1}}\neq 0\Longrightarrow M_{x,y}^{s_{1}}\neq 0.$$
Again we look at the set
$$\mathcal{E}=\{x\in\mathbb{Q}_{>0}\mid x=\pm j/i\text{ where } j<0,\ i\neq 0, Q^{i}q^{j}\in \Gamma_{+}^{s}(I) \}.$$
The largest element of $\mathcal{E}$ below $2$ is $3/2$. \\
Consider the total order given by
$$\Gamma_{+}=\{Q^{i}q^{j}\mid 3i+2j>0\} \cup \{ Q^{2j}q^{-3j}\mid j>0\}.$$
We find
$$\Gamma_{+}^{s}(I)\subseteq \{ q^{j}\mid j>0\}\cup \{Q^{i}q^{j}\mid i>0,i+j\geq 0\}$$
$$\cup \{Q^{i}q^{j}\mid j>-i>0, -j/i\geq 2\}\cup \{Q^{i}q^{j}\mid -j>i>0,-j/i\leq 3/2\}.$$
As above, we conclude that $M_{x,y}^{s_{1}}\neq 0$ for any parameter $a,b$ such that $2>a/b>3/2$.\\
We look at the set $\mathcal{E}$ (defined as above), we find that the largest element of $\mathcal{E}$ below $3/2$ is $4/3$.\\
Consider the total order given by
$$\Gamma_{+}=\{Q^{i}q^{j}\mid 4i+3j>0\} \cup \{ Q^{3j}q^{-4j}\mid j>0\}.$$
We find
$$\Gamma_{+}^{s}(I)\subseteq \{ q^{j}\mid j>0\}\cup \{Q^{i}q^{j}\mid i>0,i+j\geq 0\}$$
$$\cup \{Q^{i}q^{j}\mid j>-i>0, -j/i\geq 3/2\}\cup \{Q^{i}q^{j}\mid -j>i>0,-j/i\leq 4/3\}$$
and $M_{x,y}^{s_{1}}\neq 0$ for any parameter $a,b$ such that $3/2>a/b>4/3$.\\
We now continue the procedure. This leads us to consider the total order given by
$$\Gamma_{+}=\{Q^{i}q^{j}\mid 5i+4j>0\} \cup \{ Q^{4j}q^{-5j}\mid j>0\}.$$
We find
$$\Gamma_{+}^{s}(I)\subseteq \{ q^{j}\mid j>0\}\cup \{Q^{i}q^{j}\mid i>0,i+j\geq 0\}$$
$$\cup \{Q^{i}q^{j}\mid j>-i>0, -j/i\geq 4/3\}\cup \{Q^{i}q^{j}\mid -j>i>0,-j/i\leq 5/4\}$$
and $M_{x,y}^{s_{1}}\neq 0$ for any parameter $a,b$ such that $4/3>a/b>5/4$.\\
Finally, consider the total order given by
$$\Gamma_{+}=\{Q^{i}q^{j}\mid i+j>0\} \cup \{ Q^{j}q^{-j}\mid j>0\}.$$
We find
$$\Gamma_{+}^{s}(I)\subseteq \{ q^{j}\mid j>0\}\cup \{Q^{i}q^{j}\mid i>0,i+j\geq 0\}$$
$$\cup \{Q^{i}q^{j}\mid j>-i>0, -j/i\geq 5/4\}$$
and $M_{x,y}^{s_{1}}\neq 0$ for any parameters $a,b$ such that $5/4>a/b>1$.\\
We now calculate $M_{x,y}^{s_{1}}$ for the parameters $a,b$ where $a/b\in\{3,2, 3/2,4/3,5/4,1\}$, and we find that these are non-zero. \\
We have treated all the cases where $a\geq b$. \\
We proceed in the same way for the case $a\leq b$.\\

Doing the same for the other intervals, we can show that for $r\geq 6$ and for any parameters the coefficients
$$M_{x_{1,i}^{r},y_{1,i}^{r}}^{s_{1}}\ \text{ and }\ M_{x_{2,i}^{r},y_{2,i}^{r}}^{s_{1}}$$
are non-zero, which in turn implies that the following sets are included in a  left cell:
$$C_{i}=\{z.u_{1}^{r}.w_{i}\mid r\geq 7,z\in\Pi\},\text{ $1\leq i\leq 6$}.$$
\\
Now, let $u=s_{2}s_{1}s_{2}s_{1}s_{2}s_{3}$ and
$$W_{2}=\{e,s_{2}, s_{2}s_{1}, s_{2}s_{1}s_{2},s_{2}s_{1}s_{2}s_{1},s_{2}s_{1}s_{2}s_{1}s_{2}\}=\{w_{1},w_{2},w_{3},w_{4},w_{5},w_{6}\}.$$
For $w_{i}\in W_{2}$ let $\sigma_{w_{i}}\in \Omega$ such that $w_{i}A_{0}=A_{0}\sigma_{w_{i}}$. One can check that
$$Orb_{\Omega}(\vec{u_{1}})=\{\vec{u}_{1}\sigma_{w_{1}},...,\vec{u}\sigma_{w_{6}}\}=\{\vec{u}_{1},...,\vec{u}_{6}\}.$$
For $1\leq i\leq 6$ let
$$I_{1,i}^{r}=[s_{2}s_{3}.u^{r}.w_{i},y.u_{r}.w_{i}]=[s_{2}s_{3}.w_{i}.u_{i}^{r},y.w_{i}.u_{i}^{r}]=[x_{1,i}^{r},y_{1,i}^{r}],$$
\begin{align*}
I_{2,i}^{r}&=[s_{3}s_{2}s_{1}s_{2}s_{3}u^{r}.w_{i},s_{3}s_{2}s_{1}s_{2}s_{1}s_{2}s_{3}u^{r+1}.w_{i}]\\
&=[s_{3}s_{2}s_{1}s_{2}s_{3}.w_{i}.u_{i}^{r},s_{3}s_{2}s_{1}s_{2}s_{1}s_{2}s_{3}.w_{i}u_{i}^{r+1}]\\
&=[x_{2,i}^{r},y_{2,i}^{r}].
\end{align*}
Let $1\leq i\leq 6$ and $k=1,2$. Using our GAP3 program one can see that, for $r\geq 6,$ the intervals $I_{k,i}^{r}$  and $I_{k,i}^{r+n}$ are isomorphic for all $n\geq 0$ and that the corresponding $\mathbf{r}$-polynomials are equal. Now, using Lemma \ref{calcul}, we can show that for $r\geq 6$ and for any parameter $a,b$ such that $a/b\leq 2$ the polynomials
$$M^{s_{2}}_{x_{1,i}^{r},y_{1,i}^{r}}\ \text{ and }\ M^{s_{2}}_{x_{2,i}^{r},y_{2,i}^{r}}$$
are non-zero.\\
Therefore, for these parameters, the following sets are included in a left cell:
$$B_{i}=\{z.u^{r}.w_{i}\mid r\geq 7, z\in\Pi\},\text{ $1\leq i\leq 6$}.$$
$\ $\\

We now have a look to the parameters $a,b$ such that $a/b>2$.\\
Arguing as before, we find that, for $r\geq 6$, the polynomials
$$M_{s_{1}s_{2}s_{3}u^{r}w_{i},s_{2}s_{3}s_{2}s_{1}s_{2}s_{1}s_{2}s_{3}u^{r}w_{i}}^{s_{1}}, \ \ 
M_{s_{2}s_{3}u^{r}w_{i},s_{3}u^{r+1}w_{i}}^{s_{2}} $$
$$M_{s_{1}s_{2}s_{3}u^{r}w_{i},s_{3}s_{2}s_{1}s_{2}s_{3}u^{r}w_{i}}^{s_{1}}$$ 
are non-zero.
Therefore, for these parameters, the following sets are included in a left cell:
$$B_{i}=\{z.u^{r}.w_{i}\mid r\geq 7, z\in\Pi\},\text{ $1\leq i\leq 6$ }$$
$\ $\\

Recall that in \cite{Bremke} and \cite{Xi} it is shown that the set
$$W_{T}=\{w\in W\mid w=z'.w_{0}.z, z'\in W\}$$
is a two-sided cell and that it contains at most $\mid W_{0}\mid=12$ left cells.\\
Now one can check that the set $W_{T}$ together with the $B_{i}$'s and the $C_{i}$'s contain all the elements of $W$ except for a finite number. The theorem is proved.

Of course, we would like to find the exact decomposition of $W$ in left cells for any parameters. However, it is difficult to separate left cells. Computing some more coefficients in the case where $a>>b$, we find a more precise decomposition of $W$ which is included in the left cells decomposition. 

We show this decomposition in Figure 2. We identify $w\in W$ with the alcove $wA_{0}$. The sets which are included in a left cells are formed by the alcoves lying in the same connected component after removing the thick line.\\
We have
$$W_{T}=\overset{12}{\underset{i=1}{\cup}}A_{i}.$$
The figure also show the shape of the sets $B_{i}$ and $C_{i}$. 
\begin{center}
\begin{pspicture}(-6,6)(8.2,-8.2)
\psset{unit=0.85cm}
\SpecialCoor
\psline(0,-6)(0,6)
\multido{\n=1+1}{7}{
\psline(\n;30)(\n;330)}
\rput(1;30){\psline(0,0)(0,5.5)}
\rput(1;330){\psline(0,0)(0,-5.5)}
\rput(2;30){\psline(0,0)(0,5)}
\rput(2;330){\psline(0,0)(0,-5)}
\rput(3;30){\psline(0,0)(0,4.5)}
\rput(3;330){\psline(0,0)(0,-4.5)}
\rput(4;30){\psline(0,0)(0,4)}
\rput(4;330){\psline(0,0)(0,-4)}
\rput(5;30){\psline(0,0)(0,3.5)}
\rput(5;330){\psline(0,0)(0,-3.5)}
\rput(6;30){\psline(0,0)(0,3)}
\rput(6;330){\psline(0,0)(0,-3)}
\rput(7;30){\psline(0,0)(0,2.5)}
\rput(7;330){\psline(0,0)(0,-2.5)}
\multido{\n=1+1}{7}{
\psline(\n;150)(\n;210)}
\rput(1;150){\psline(0,0)(0,5.5)}
\rput(1;210){\psline(0,0)(0,-5.5)}
\rput(2;150){\psline(0,0)(0,5)}
\rput(2;210){\psline(0,0)(0,-5)}
\rput(3;150){\psline(0,0)(0,4.5)}
\rput(3;210){\psline(0,0)(0,-4.5)}
\rput(4;150){\psline(0,0)(0,4)}
\rput(4;210){\psline(0,0)(0,-4)}
\rput(5;150){\psline(0,0)(0,3.5)}
\rput(5;210){\psline(0,0)(0,-3.5)}
\rput(6;150){\psline(0,0)(0,3)}
\rput(6;210){\psline(0,0)(0,-3)}
\rput(7;150){\psline(0,0)(0,2.5)}
\rput(7;210){\psline(0,0)(0,-2.5)}
\multido{\n=1.5+1.5}{4}{
\psline(-6.062,\n)(6.062,\n)}
\multido{\n=0+1.5}{5}{
\psline(-6.062,-\n)(6.062,-\n)}
\psline(0;0)(7;30)

\rput(0,1){\psline(0;0)(7;30)}
\rput(0,1){\psline(0;0)(7;210)}
\rput(0,2){\psline(0;0)(7;30)}
\rput(0,2){\psline(0;0)(7;210)}
\rput(0,3){\psline(0;0)(6;30)}
\rput(0,3){\psline(0;0)(7;210)}
\rput(0,4){\psline(0;0)(4;30)}
\rput(0,4){\psline(0;0)(7;210)}
\rput(0,5){\psline(0;0)(2;30)}
\rput(0,5){\psline(0;0)(7;210)}
\rput(0,6){\psline(0;0)(7;210)}
\rput(-1.732,6){\psline(0;0)(5;210)}
\rput(-3.464,6){\psline(0;0)(3;210)}
\rput(-5.196,6){\psline(0;0)(1;210)}

\rput(0,-1){\psline(0;0)(7;30)}
\rput(0,-1){\psline(0;0)(7;210)}
\rput(0,-2){\psline(0;0)(7;30)}
\rput(0,-2){\psline(0;0)(7;210)}
\rput(0,-3){\psline(0;0)(7;30)}
\rput(0,-3){\psline(0;0)(6;210)}
\rput(0,-4){\psline(0;0)(7;30)}
\rput(0,-4){\psline(0;0)(4;210)}
\rput(0,-5){\psline(0;0)(7;30)}
\rput(0,-5){\psline(0;0)(2;210)}
\rput(0,-6){\psline(0;0)(7;30)}
\rput(1.732,-6){\psline(0;0)(5;30)}
\rput(3.464,-6){\psline(0;0)(3;30)}
\rput(5.196,-6){\psline(0;0)(1;30)}

\psline(0;0)(6.928;60)
\rput(1.732,0){\psline(0;0)(6.928;60)}
\rput(3.464,0){\psline(0;0)(5.196;60)}
\rput(5.196,0){\psline(0;0)(1.732;60)}
\rput(1.732,0){\psline(0;0)(6.928;240)}
\rput(3.464,0){\psline(0;0)(6.928;240)}
\rput(5.196,0){\psline(0;0)(6.928;240)}
\rput(6.062,-1.5){\psline(0;0)(5.196;240)}
\rput(6.062,-4.5){\psline(0;0)(1.732;240)}
\rput(-1.732,0){\psline(0;0)(6.928;60)}
\rput(-3.464,0){\psline(0;0)(6.928;60)}
\rput(-5.196,0){\psline(0;0)(6.928;60)}
\rput(-1.732,0){\psline(0;0)(6.928;240)}
\rput(-3.464,0){\psline(0;0)(5.196;240)}
\rput(-5.196,0){\psline(0;0)(1.732;240)}
\rput(-6.062,1.5){\psline(0;0)(5.196;60)}
\rput(-6.062,4.5){\psline(0;0)(1.732;60)}

\psline(0;0)(6.928;120)
\psline(0;0)(6.928;300)
\rput(-1.732,0){\psline(0;0)(6.928;120)}
\rput(-3.464,0){\psline(0;0)(5.196;120)}
\rput(-5.196,0){\psline(0;0)(1.732;120)}
\rput(-1.732,0){\psline(0;0)(6.928;300)}
\rput(-3.464,0){\psline(0;0)(6.928;300)}
\rput(-5.196,0){\psline(0;0)(6.928;300)}
\rput(-6.062,-1.5){\psline(0;0)(5.196;300)}
\rput(-6.062,-4.5){\psline(0;0)(1.732;300)}
\rput(1.732,0){\psline(0;0)(6.928;300)}
\rput(3.464,0){\psline(0;0)(5.196;300)}
\rput(5.196,0){\psline(0;0)(1.732;300)}
\rput(1.732,0){\psline(0;0)(6.928;120)}
\rput(3.464,0){\psline(0;0)(6.928;120)}
\rput(5.196,0){\psline(0;0)(6.928;120)}
\rput(6.062,1.5){\psline(0;0)(5.196;120)}
\rput(6.062,4.5){\psline(0;0)(1.732;120)}

\rput(0,1){\psline(0;0)(7;150)}
\rput(0,1){\psline(0;0)(7;330)}
\rput(0,2){\psline(0;0)(7;150)}
\rput(0,2){\psline(0;0)(7;330)}
\rput(0,3){\psline(0;0)(6;150)}
\rput(0,3){\psline(0;0)(7;330)}
\rput(0,4){\psline(0;0)(4;150)}
\rput(0,4){\psline(0;0)(7;330)}
\rput(0,5){\psline(0;0)(7;330)}
\rput(0,6){\psline(0;0)(7;330)}
\rput(-1.732,6){\psline(0;0)(2;330)}

\rput(0,-1){\psline(0;0)(7;150)}
\rput(0,-1){\psline(0;0)(7;330)}
\rput(0,-2){\psline(0;0)(7;150)}
\rput(0,-2){\psline(0;0)(7;330)}
\rput(0,-3){\psline(0;0)(7;150)}
\rput(0,-3){\psline(0;0)(6;330)}
\rput(0,-4){\psline(0;0)(7;150)}
\rput(0,-4){\psline(0;0)(4;330)}
\rput(0,-5){\psline(0;0)(7;150)}
\rput(0,-5){\psline(0;0)(2;330)}
\rput(0,-6){\psline(0;0)(7;150)}

\rput(6.062,3.5){\psline(0;0)(5;150)}
\rput(6.062,4.5){\psline(0;0)(3;150)}
\rput(6.062,5.5){\psline(0;0)(1;150)}

\rput(-6.062,-3.5){\psline(0;0)(5;330)}
\rput(-6.062,-4.5){\psline(0;0)(3;330)}
\rput(-6.062,-5.5){\psline(0;0)(1;330)}

\psline(0;0)(7;330)
\psline(0;0)(7;150)
\psline(0;0)(7;210)
\psline(0;0)(6.928;240)

\pspolygon[linewidth=0.7mm](0,0)(0,1)(0.433,0.75)
\pspolygon[linewidth=0.7mm](0,0)(0,1)(-0.866,1.5)
\pspolygon[linewidth=0.7mm](0,0)(-0.433,0.75)(-1.732,0)(-0.866,0)(-0.866,-1.5)
\pspolygon[linewidth=0.7mm](0,0)(0.433,0.75)(1.732,0)(0.866,0)(0.866,-1.5)
\pspolygon[linewidth=0.7mm](0,1)(0,1.5)(-0.866,1.5)
\pspolygon[linewidth=0.7mm](0,1)(0,3)(0.433,2.25)(1.732,3)(0.433,0.75)
\psline[linewidth=0.7mm](0,0)(0,-6)
\psline[linewidth=0.7mm](0;0)(6.928;240)
\rput(0.866,-1.5){\psline[linewidth=0.7mm](0;0)(4.5;270)}
\rput(0.866,-1.5){\psline[linewidth=0.7mm](0;0)(5.196;300)}
\rput(1.732,0){\psline[linewidth=0.7mm](0;0)(6.928;300)}
\rput(2.598,-1.5){\psline[linewidth=0.7mm](0;0)(4;330)}
\rput(1.732,0){\psline[linewidth=0.7mm](0;0)(5;330)}
\rput(1.732,0){\psline[linewidth=0.7mm](0;0)(4.33;0)}
\rput(0.866,1.5){\psline[linewidth=0.7mm](0;0)(5.196;0)}
\rput(0.866,1.5){\psline[linewidth=0.7mm](0;0)(6;30)}
\rput(2.598,1.5){\psline[linewidth=0.7mm](0;0)(4;30)}
\rput(1.732,3){\psline[linewidth=0.7mm](0;0)(3.464;60)}
\rput(0,3){\psline[linewidth=0.7mm](0;0)(3.464;60)}
\rput(0.866,4.5){\psline[linewidth=0.7mm](0;0)(1.5;90)}
\rput(0,3){\psline[linewidth=0.7mm](0;0)(3;90)}
\rput(0,3){\psline[linewidth=0.7mm](0;0)(3.464;120)}
\rput(-0.866,1.5){\psline[linewidth=0.7mm](0;0)(5.196;120)}
\rput(-0.866,1.5){\psline[linewidth=0.7mm](0;0)(6;150)}
\rput(-0.866,1.5){\psline[linewidth=0.7mm](0;0)(5.196;180)}
\rput(-2.598,1.5){\psline[linewidth=0.7mm](0;0)(4;150)}
\rput(-1.732,0){\psline[linewidth=0.7mm](0;0)(4.33;180)}
\rput(-1.732,0){\psline[linewidth=0.7mm](0;0)(5;210)}
\rput(-1.732,0){\psline[linewidth=0.7mm](0;0)(6.928;240)}
\rput(-1.732,0){\psline[linewidth=0.7mm](0;0)(6.928;240)}
\rput(-2.598,-1.5){\psline[linewidth=0.7mm](0;0)(4;210)}
\rput(0,1.5){\psline[linewidth=0.7mm](0;1.5)(0.866;1.5)}

\rput(0.17,0.68){{\tiny $A_{0}$}}
\rput(-1.5,-4.33){{\small $A_{1}$}}
\rput(-4.97,-4.3){{\small $A_{2}$}}
\rput(-4.97,-1.3){{\small $A_{3}$}}
\rput(-4.97,1.7){{\small $A_{4}$}}
\rput(-2.372,3.2){{\small $A_{5}$}}
\rput(-0.23,4.67){{\small $A_{6}$}}
\rput(1.1,5.85){{\small $A_{7}$}}
\rput(2.824,3.2){{\small $A_{8}$}}
\rput(3.69,1.7){{\small $A_{9}$}}
\rput(5.452,-1.31){{\small $A_{10}$}}
\rput(5.452,-4.31){{\small $A_{11}$}}
\rput(1.988,-4.31){{\small $A_{12}$}}
\psline[linewidth=0.5mm]{->}(0.433,5.75)(0.433,6.75)
\rput(0.433,7){$C_{1}$}
\rput(5.812,3.8){\psline[linewidth=0.5mm]{->}(0;0)(1;30)}
\rput(7,4.4){$C_{2}$}
\rput(5.812,-2.8){\psline[linewidth=0.5mm]{->}(0;0)(1;330)}
\rput(7,-3.4){$C_{4}$}
\psline[linewidth=0.5mm]{->}(0.433,-5.75)(0.433,-6.75)
\rput(0.433,-7){$C_{6}$}
\rput(-5.812,-2.8){\psline[linewidth=0.5mm]{->}(0;0)(1;210)}
\rput(-7,-3.4){$C_{5}$}
\rput(-5.812,3.8){\psline[linewidth=0.5mm]{->}(0;0)(1;150)}
\rput(-7,4.4){$C_{3}$}
\rput(2.4,5.75){\psline[linewidth=0.5mm]{->}(0;0)(1;60)}
\rput(3.05,6.83){$B_{1}$}
\rput(5.812,0.75){\psline[linewidth=0.5mm]{->}(0;0)(1;0)}
\rput(7.1,0.75){$B_{3}$}
\rput(4.132,-5.75){\psline[linewidth=0.5mm]{->}(0;0)(1;300)}
\rput(4.837,-6.83){$B_{5}$}
\rput(-4.132,-5.75){\psline[linewidth=0.5mm]{->}(0;0)(1;240)}
\rput(-4.837,-6.83){$B_{6}$}
\rput(-5.812,0.75){\psline[linewidth=0.5mm]{->}(0;0)(1;180)}
\rput(-7.1,0.75){$B_{4}$}
\rput(-2.4,5.75){\psline[linewidth=0.5mm]{->}(0;0)(1;120)}
\rput(-3.05,6.83){$B_{2}$}

\rput(0,-7.8){Fig. 2.}

\end{pspicture}
\end{center}


\begin{thebibliography}{10}

\bibitem{Bed}
R.~Bedard.
\newblock Cells for two {C}oxeter groups.
\newblock {\em Comm. Algebra $\mathbf{14}$, 1253-1286}, 1986.

\bibitem{Bremke}
K.~Bremke.
\newblock On generalized cells in affine {W}eyl groups.
\newblock {\em Journal of Algebra $\mathbf{191}$, 149-173}, 1997.

\bibitem{Chen}
C.~Chen.
\newblock The decomposition into left cells of the affine {W}eyl group of type
  $\tilde{D}_{4}$.
\newblock {\em Journal of Algebra $\mathbf{163}$, 692-728}, 1994.

\bibitem{Du}
J.~Du.
\newblock The decomposition into cells of the affine {W}eyl group of type
  $\tilde{B}_{3}$.
\newblock {\em Comm. Algebra $\mathbf{16}$, 1383-1409}, 1988.

\bibitem{Cox}
F.~Du Cloux.
\newblock Coxeter, version 3.
\newblock {\em available from \\
http://math.univ-lyon1.fr/$\sim$ducloux/coxeter/coxeter3/english/coxeter3\_e.html}


\bibitem{geck}
M.~Geck.
\newblock Computing {K}azhdan-{L}usztig cells for unequal parameters.
\newblock {\em Journal of Algebra $\mathbf{281}$, 342-365}, 2004.

\bibitem{Lus1}
G.~Lusztig.
\newblock Hecke algebras and {J}antzen's generic decomposition patterns.
\newblock {\em Advances in Mathematics $\mathbf{37}$, 121-164}, 1980.

\bibitem{Lus1p}
G.~Lusztig.
\newblock Left cells in {W}eyl groups.
\newblock {\em Lectures Notes in Math., Springer, Vol. 1024, 99-111}, 1983.

\bibitem{Lus3}
G.~Lusztig.
\newblock Cells in affine {W}eyl groups.
\newblock {\em Advanced Studies in Pure Math., Vol. 6, 255-287}, 1985.

\bibitem{Lus2}
G.~Lusztig.
\newblock The two-sided cells of the affine {W}eyl group of type
  $\tilde{A}_{n}$.
\newblock {\em Math. Sci. Res. Inst. Publ., 4, 275-283}, 1985.

\bibitem{Lus4}
G.~Lusztig.
\newblock Cells in affine {W}eyl groups {I}{I}.
\newblock {\em Journal of Algebra $\mathbf{109}$, 536-548}, 1987.

\bibitem{Lus5}
G.~Lusztig.
\newblock Cells in affine {W}eyl groups {I}{I}{I}.
\newblock {\em J. Fac. Sci. Univ. Tokyo Sect. IA, Math. $\mathbf{34}$,
  223-243}, 1987.

\bibitem{Lus6}
G.~Lusztig.
\newblock Cells in affine {W}eyl groups {I}{V}.
\newblock {\em J. Fac. Sci. Univ. Tokyo Sect. IA, Math. $\mathbf{36}$,
  297-328}, 1989.

\bibitem{bible}
G.~Lusztig.
\newblock {H}ecke algebras with unequal parameters.
\newblock {\em CRM Monographs Ser. $\mathbf{18}$, Amer. Math. Soc. ,
  Providence, RI}, 2003.

\newcommand{\etalchar}[1]{$^{#1}$}

  \bibitem{GAP}
            Martin Sch{\accent127 o}nert et~al.
  \newblock {\em {GAP} --
            {Groups}, {Algorithms}, and {Programming} --
            version 3 release 4 patchlevel 4"}.
  \newblock Lehrstuhl D f{\accent127 u}r Mathematik,
            Rheinisch Westf{\accent127 a}lische
            Technische Hoch\-schule, Aachen, Germany, 1997.   

\bibitem{Shi1}
J.-Y. Shi.
\newblock The {K}azhdan-{L}usztig cells in certain affine {W}eyl groups.
\newblock {\em Lectures Notes in Math., Springer-Verlag, Vol. 1179}, 1986.

\bibitem{Shi3}
J.-Y. Shi.
\newblock Left cells in affine {W}eyl group ${W}_{a}(\tilde{D}_{4})$.
\newblock {\em Osaka J. Math. $\mathbf{31}$, 27-50}, 1994.

\bibitem{Shi2}
J.-Y. Shi.
\newblock Left cells in affine {W}eyl groups.
\newblock {\em Tokohu Math. J. $\mathbf{46}$, 105-124}, 1994.

\bibitem{Xi}
N.~Xi.
\newblock Representations of affine {H}ecke algebras.
\newblock {\em Lectures Notes in Math., Springer-Verlag, Vol. 1587}, 1994.

\end{thebibliography}
\end{document}